\documentclass[11pt]{article}
\usepackage{amssymb}
\usepackage{amsmath}
\usepackage{amsthm}
\newcommand{\bnu}{\boldmath$\nu$}
\newcommand{\mbnu}{\mbox{\bnu}}
\newcommand{\bomega}{\boldmath$\omega$}
\newcommand{\mbomega}{\mbox{\bomega}}
\newcommand{\bDelta}{\boldmath$\Delta$}
\newcommand{\mbDelta}{\mbox{\bDelta}}
\newcommand{\bLambda}{\boldmath$\Lambda$}
\newcommand{\mbLambda}{\mbox{\bLambda}}
\newcommand{\bGamma}{\boldmath$\Gamma$}
\newcommand{\mbGamma}{\mbox{\bGamma}}
\newcommand{\bSigma}{\boldmath$\Sigma$}
\newcommand{\mbSigma}{\mbox{\bSigma}}
\newcommand{\btheta}{\boldmath$\theta$}
\newcommand{\mbtheta}{\mbox{\btheta}}

\newcommand{\bbeta}{\boldmath$\beta$}
\newcommand{\mbbeta}{\mbox{\bbeta}}
\newcommand{\bmu}{\boldmath$\mu$}
\newcommand{\mbmu}{\mbox{\bmu}}

\newcommand{\reg}{$\emptyset \subseteq a \subseteq P$}
\newcommand{\rag}{$P_{1} \subseteq a \subseteq P$}
\newcommand{\Xbar}{$\bar{\bf X}$}
\newcommand{\Ybar}{$\bar{\bf Y}$}
\newcommand{\mXbar}{\mbox{\Xbar}}
\newcommand{\mYbar}{\mbox{\Ybar}}
\newcommand{\Zbar}{$\bar{\bf Z}$}
\newcommand{\mZbar}{\mbox{\Zbar}}

\newfont{\bfit}{cmbxti10 scaled\magstep1}

\DeclareMathOperator{\tr}{tr}

\begin{document}

\newtheorem{theorem}{\indent {\sc Theorem}}
\newtheorem{prop}{\indent {\sc Proposition}}
\newtheorem{lemma}{\indent {\sc Lemma}}
\renewcommand{\proofname}{\hspace*{\parindent}{\sc Proof.}}
\thispagestyle{empty}

\begin{center}
{\bf TOWARDS AN UNIFIED THEORY FOR TESTING STATISTICAL HYPOTHESES: MULTIVARIATE MEAN WITH NUISANCE
COVARIANCE MATRIX}\\
\medskip
{\ B{\sc y}  M{\sc ing}-T{\sc ien} T{\sc sai}} \\
\medskip
{\it Academia Sinica}
\end{center}

\vspace{0.3cm}
\begin{center}
\small
\parbox{15cm}{\sloppy \ \,
Under a multinormal distribution with arbitrary unknown covariance matrix, the main purpose of this paper is to
propose a framework to achieve the goal of reconciliation of Bayesian, frequentist and Fisherian paradigms for
the problems of testing mean against restricted alternatives (closed convex cones). Combine fiducial inference
and Wald decision theory via d-admissibility into an unified approach, the goal can then be achieved. To proceed,
the tests constructed via the union-intersection principle (Roy, 1953) are studied. The difficulty for likelihood
ratio principle is mentioned. }

\end{center}

\vspace{1cm}
\noindent \underline{~~~~~~~~~~~~~~~~~~~~~~~~~~~~~~~~~~~~~~~~~~~~~~~~~~~~~~~~~}\\
\indent{\small{\it AMS 1991 subject classification}: 62C15, 62H15.}

\newpage

\vspace{0.3cm}
\indent {\bf 1. Introduction and preliminary notions.} Let \{${\bf X}_{i}; 1 \leq i \leq n$\}
be independent and identically distributed random vectors (i.i.d.r.v.) having a $p$-variate
normal distribution with mean vector {\bmu} and dispersion matrix {\mbSigma}. Consider the
null hypothesis
\def \theequation{1.\arabic{equation}}
\setcounter{equation}{0}
\begin{eqnarray}
H_{0}^{\ast}:{\mbmu} \in \Gamma_{1} = \{
{\mbmu}|{\bf B}_{1}{\mbmu} = {\bf 0}\}
\end{eqnarray}
against a restricted (convex polyhedral cone) alternative
\begin{eqnarray}
H_{1}^{\ast}:{\mbmu} \in \Gamma_{2} = \{
{\mbmu}|{\bf B}_{2}{\mbmu} \geq {\bf 0}\},
\end{eqnarray}
where ${\bf B}_{i} \in {\cal C}(m_{i},p)$, the set of $m_{i}\times p$ matrices of
rank $m_{i} (1 \leq m_{i}\leq p),i$=1,2, and $\Gamma_{1}$ is
assumed to be the linear hull of ${\tilde{\Gamma}}_{2} = \{
{\mbmu}|{\bf B}_{1}{\mbmu} = {\bf 0},
{\bf B}_{2}{\mbmu} \geq {\bf 0}\}$. An impasse in a general formulation of optimal tests for
this problem is the lack of invariance of the model (as well as the likelihood ratio or
allied tests) under suitable groups of transformations which map the sample space onto
itself. Specificially, if we consider an affine transformation ${\bf X} \rightarrow {\bf Y}=
{\bf a} + {\bf AX}$, where ${\bf A}$ is non-singular, this hypothesis testing problem is not
invariant in general although it is invariant under two special groups of linear
transformations (positive diagonal matrices and permutation matrices) when
${\bf B}_{1}={\bf B}_{2}={\bf I}$. As such, a canonical reduction of the noncentrality to a
single coordinate may not work out, and if lacking this invariance, the usual techniques fail
to provide an optimality property of the tests. For this reason, often, a hypothesis
related transformation is used on the basic r.v.'s, and on this transformed r.v.'s suitable
tests are formulated. Towards this, we may consider
${\bf Y}_{i} = {\bf B}_{1}{\bf X}_{i},i \geq 1,{\mbtheta} =
{\bf B}_{1}{\mbmu}$ and ${\bf B} = {\bf B}_{2}{\bf B}_{1}^{'}{({\bf B}_{1}{\bf B}_
{1}^{'})}^{-1}$, so that based on the ${\bf Y}_{i}$, the hypotheses in (1.1)
and (1.2) can be expressed as
\begin{eqnarray}
{\bar{H}}_{0}:{\mbtheta} = {\bf 0}~~\hbox{vs.}~~{\bar{H}}_{1}:{\mbtheta} \in
\Gamma_{3} = \{{\mbtheta} | {\bf B}{\mbtheta} \geq {\bf 0}, ||
{\mbtheta}|| > 0\},
\end{eqnarray}
(Recall that ${\mbSigma} \to {\mbSigma}_{1}={\bf B}_{1}{\mbSigma}{\bf B}_{1}^{'}$ and need
not be diagonal). Further reduction to a positive orthant model is feasible by the
transformations ${\bf Z}_{i}={\bf B}{\bf Y}_{i}, {\mbbeta}={\bf B}{\mbtheta}$ and
${\mbSigma} \to {\bf B}{\mbSigma}_{1}{\bf B}^{'}$ . Situations are similar for the cases
when the conditions $m_{i}\leq p,i=1,2$ are relaxed. Also note that a specific halfspace
can be transformed into another halfspace by a non-singular linear transformation. Hence
without loss of generality it suffices to consider the following models:

\indent Let ${\bf X}_{1},\ldots,{\bf X}_{n}$ be i.i.d.r.v.'s
with the $N_{p}({\mbtheta},{\mbSigma})$ d.f., where ${\mbSigma}$ is unknown and assumed to
be positive definite (p.d.). Consider the hypotheses
\begin{eqnarray}
H_{0}:{\mbtheta} = {\bf 0}~~\hbox{vs.}~~H_{1}:{\mbtheta} \in {\cal C} \backslash \{{\bf 0}\},
\end{eqnarray}
where ${\cal C}$ denotes a closed convex cone containing a p-dimensional open set. Denote
the positive orthant space by
${\cal O}_{p}^{+} = \{{\mbtheta}\in R^{p}|~{\mbtheta} \geq {\bf 0}\}$.
Notice that when ${\cal C}$ is a proper set contained in a halfspace, under a suitable
linear transformation the problem in (1.4) can be reduced to the problem for testing against
the positive orthant space with another unknown positive definite covariance matrix. When
${\cal C}$ is a specific halfspace, then it can be transformed into another halfspace by a
non-singular linear transformation. Hence without loss of generality it suffices enough to
study the cases that ${\cal C}$ is the positive orthant space ${\cal O}_{p}^{+}$ and
${\cal C}$ is the halfspace ${\cal H}_{p}^{*}=\{{\mbtheta}\in R^{p}|~{\theta}_{p} \geq 0 \}$
throughout this paper.

\indent For testing against global alternative,
$H^{g}_{0}:{\mbtheta} = {\bf 0}~~\hbox{vs.}~~H^{g}_{1}:{\mbtheta}
\ne {\bf 0}$, the likelihood ratio test (LRT) and the union-intersection test (UIT) are isomorphic, and are
well-known as the Hotelling's $T^{2}$-test (Anderson, 2003). Kiefer
and Schwartz (1965) showed that the Hotelling's $T^{2}$-test is a
proper Bayes test. The Hotelling's $T^{2}$-test is also known
as a version of integrated LRT with respect to objective priors. As
such, for testing against global alternative, there are several ways
to establish the equivalent relationships of the Hotelling's
$T^{2}$-test with Bayes tests.

\indent On the other hand, for the problems of testing against restricted alternatives
considered in (1.4), the LRT and the UIT are different (Perlman, 1969, Sen and Tsai, 1999)
when the covaraince matrix ${\mbSigma}$ is totally unknown. The problem (1.4) indeed provides
us some partial informations of interesting parameters. However, due to the difficulty of
integration over high-dimensional restricted parameter space the explicit forms of Bayes tests
are hardly obtained, though the numerical values of them can be obtained via the method of
Markov chain Monte Carlo. The same difficulty arises for the approach of Bayes factor.
In Section 2, we show that the UIT are no longer the Bayes tests for the problem (1.4).
The likelihood integrated (Berger, Liseo and Wolpert, 1999) with respect to a relevant Haar
measure for the nuisance covariance matrix while remaining neutral with respect to the
parameter of interest, which fiducial inference was intended to be, can be viewed in a objective
Bayesian light. It is shown that for the problems of testing against restricted alternatives
both the LRT and the UIT are the versions of integrated LRT and integrated UIT, respectively.
These result in the reconciliation of frequentist and Bayesian paradigms.

\indent For testing against global alternative, the Hotelling's $T^{2}$-test enjoys many optimal
properties of Neyman-Pearson testing hypothesis theory such as similarity, unbiasedness,
power monotonicity, most stringent, uniformly most powerful invariant and $\alpha$-admissibility
(Anderson, 2003). For the problem of testing against the positive orthant space, both the LRT
and the UIT are shown not unbiased (Perlman, 1969, Sen and Tsai, 1999). And it is easy to
see that both the LRT and the UIT are power dominated by the corresponding LRT and UIT for the problem
of testing against the halfsapce ${\cal H}^{*}_{p}$, respectively. As such, for the problem of
testing against the postive orthant space both the LRT and the UIT are $\alpha$-inadmissible.
The related domination problems are also studied in Section 3, some of these domination results
are against the common statistical sense. The situation in Neyman-Pearson theory is beautifully
peaceful and compelling for the unrestricted alternative, while it is not in some restricted
alternatives world. One of the major reasons for these phenomena is that for the problem of
testing against the positive orthant space the distribution functions of both the LRT statistic
and the UIT statistic are not free from the dependence on the unkown nuisance covariance matrix
under null hypothesis. Because of this unpleasent feature, the Fisher's approach of reporting
$p$-values may also not work well for the problem of testing against the positive orthant space.
Several methods have been proposed to get rid of the unpleasent feature, we mention some in
Section 3.

\indent In the literature, there are many ways to eliminate the nuisance parameters (Basu, 1977).
Most of them can be applied to the problem of testing against the positive orthant space. As we
know from Section 3 that the domination problems of hypothesis testing heavily depend on the
choice of the critical points which are used to report the $p$-values. Sometimes the danger
stems from too narrow definitions of what is meant by optimality and the strongest intuition
can sometimes go astray too. Hence, the question naturally raised is whether there exists any
satisfactory, unified approach to overcome the difficulties. Neither Fisher's approach of
reporting $p$-values alone nor Neyman-Pearson's optimal theory for power function alone is a
well statisfactory criterion for evaluating the performance of tests. The spirit of compromise
between Fisher's approach and Neyman-Pearson's optimal theory without detailed consideration of
power may shed light on the testing hypothesis theory. Imposing on the balance between type I
error and power, Wald's decision theory still paves a unified way to combine the best features
of both Neyman-Pearson's and Fisher's ideas. In Section 4, we show that the Hotelling's
$T^{2}$-test is inadmissible for the problem (1.4). And further show that the UIT is
$\alpha$-admissible and $d$-admissible for the problem of testing against the halfspace, and
it is $\alpha$-inadmissible but $d$-admissible for the problem of testing against the positive
orthant space, respectively. We fail to claim the same results for LRT in the problem (1.4)
due to the facts that the acceptance region of LRT is not convex.

\vspace{0.3cm}
\def \theequation{2.\arabic{equation}}
\setcounter{equation}{0}

\indent {\bf 2. A reconciliation of frequentist and Bayesian paradigms.} For every $n$ $(\geq 2)$, let
\begin{eqnarray}
  {\mXbar}_{n} = n^{-1} \sum_{i=1}^{n} {\bf X}_{i}~~\hbox{and}~~{\bf S}_{n} =(n-1)^{-1}
  \sum_{i=1}^{n} ({\bf X}_{i} - {\mXbar}_{n}) ({\bf X}_{i} -
  {\mXbar}_{n})',
\end{eqnarray}
the Hotelling's $T^{2}$-test is expressed as
\begin{eqnarray}
T^{2}_{n}=n{\mXbar}'_{n}{\bf S}^{-1}_{n}{\mXbar}_{n}.
\end{eqnarray}
For testing against global alternative, Kiefer and Schwartz (1965) showed that the Hotelling's
$T^{2}$-test is a proper Bayes test. The main goal of this section is
to see whether the reconciliation of frequentist and Bayesian paradigms can be established for
the problem (1.4). First, we establish the necessary conditions for the Bayes tests.

\vspace{0.3cm}
\indent {\sc Proposition 2.1}. {\it Let ${\mathcal A}$ denote the acceptance region of a
size-$\alpha$ Bayes test for the problem of testing $H_{0}: {\mbtheta} = {\bf 0}$ against
$H_{1}: {\mbtheta} \in {\cal C} \backslash \{{\bf 0}\}$, where ${\cal C}$ is either the
positive orthant space ${\cal O}_{p}^{+}$ or the halfspace ${\cal H}_{p}^{*}$, and
${\mathcal L}$ be any line of support of ${\mathcal A}$. Denote
$\Bar{\mathcal A}={\mathcal A}\cup \partial {\mathcal A}$, $\partial {\mathcal A}$ being the
boundary set of ${\mathcal A}$. Then either ${\mathcal L} \subseteq \partial {\mathcal A}$
or ${\mathcal L} \cap \Bar{\mathcal A}=\{{\bf a}\}$, where
${\bf a} \in \partial {\mathcal A}$.}

\indent {\sc Proof}. The density of ${\bf X}_{1}, \ldots, {\bf X}_{n}$ is
\begin{eqnarray}
 \frac{e^{\frac{-1}{2} n{\mbtheta}^{'}{\mbSigma}^{-1}{\mbtheta}}}
   {(2\pi)^{\frac{1}{2}pn}|{\mbSigma}|^{\frac{1}{2}n}}
   \mbox{exp}\left[{n}{\mbtheta}^{'}{\mbSigma}^{-1}{\mXbar}_{n}
    +\mbox{tr}(-\frac{1}{2}{\mbSigma}^{-1})
    \sum_{j=1}^{n}{\bf X}_{j}{\bf X}^{'}_{j}\right].
\end{eqnarray}
The vector ${\bf y}=({\bf y}^{(1)'},{\bf y}^{(2)'})'$ is composed of
${\bf y}^{(1)}={\mXbar}_{n}$ and ${\bf y}^{(2)}=(v_{11}, 2v_{12},\cdots,2v_{1p},
  v_{22},\cdots, v_{pp})$,
where $\left(v_{ij}\right)={\bf V}_{n}=\sum_{j=1}^{n}{\bf X}_{j}{\bf X}^{'}_{j}$.
The vector ${\mbomega}=({\mbomega}^{(1)'},{\mbomega}^{(2)'})'$ is composed of
${\mbomega}^{(1)}=n{\mbSigma}^{-1}{\mbtheta}$ and ${\mbomega}^{(2)}
=-\frac{1}{2}({\sigma}^{11}, {\sigma}^{12},\cdots,{\sigma}^{1p}, {\sigma}^{22}, \cdots,
{\sigma}^{pp})'$, where $({\sigma}^{ij})={\mbSigma}^{-1}$. Write
${\mbSigma}^{-1}={\mbGamma}({\mbomega}^{(2)})$, thus the density of ${\bf Y}$ becomes
\begin{eqnarray}
\frac{|{\mbGamma}({\mbomega}^{(2)})|^{\frac{1}{2}n}\mbox{e}^{-\frac{1}{2}
{\mbomega}^{(1)'}{\mbomega}^{(1)}}}{(2\pi)^{\frac{1}{2}pn}}
  \mbox{exp} \{{\mbomega}^{'}{\bf y}\}
\end{eqnarray}
Since $\mbSigma$ is arbitrary positive definite, ${\mbtheta}\in {\cal O}_{p}^{+}$ implies
that ${\mbomega}^{(1)} \in {\cal H}^{+}_{p}$, where
${\cal H}^{+}_{p}=\{{\mbtheta}\in R^{p}|~\sum_{i=1}^{p}{\theta}_{i} \geq 0 \}$.
Similarly, let ${\cal S}$ be the set so that
${\mbomega}^{(1)} \in {\cal S}$ whenever ${\mbtheta}\in {\cal H}_{p}^{*}$. Thus for
simplicity and without loss of generality, it suffices to work the case that
${\mbomega}^{(1)} \in {\cal H}^{+}_{p}$. Let
${\mathcal M}$ be the group of positive definite matrices and $G(\mbomega)$ denote an
arbitrary prior distribution function $G(\mbomega)$ on the set
${\mathcal N}=\{\mbomega|~{\mbomega}^{(1)}\in~{\cal H}^{+}_{p}~\hbox{and}~
{\mbGamma}({\mbomega}^{(2)}) \in {\mathcal M}\}$. Let $G_0(\mbomega)$ denote the prior
distribution function on the set
${\mathcal N}_0=\{\mbomega|~{\mbomega}^{(1)}={\bf 0}~\hbox{and}~
{\mbGamma}({\mbomega}^{(2)}) \in {\mathcal M}\}$. Then the size-$\alpha$ Bayes test of
$H_0$ vs. $H_1$ against this prior distribution rejects $H_0$ if
\begin{eqnarray}
t({\bf y})=\int_{{\mathcal N}}\mbox{exp}\{{\mbomega}^{(1)'}{\bf y}^{(1)}\}
   dG^{*}({\mbomega}) \geq c,
\end{eqnarray}
where
\begin{eqnarray}
dG^{*}({\mbomega})=\frac{\mbox{exp}\{{\mbomega}^{(2)'}{\bf y}^{(2)}\}|
   {\mbGamma}({\mbomega}^{(2)})|^{\frac{1}{2}n}\mbox{e}^{-\frac{1}{2}{\mbomega}^{(1)'}
   {\mbomega}^{(1)}}dG({\mbomega})}{\int_{{\mathcal N}_0}
   \mbox{exp}\{{\mbomega}^{(2)'}{\bf y}^{(2)}\}|{\mbGamma}
   ({\mbomega}^{(2)})|^{\frac{1}{2}n}dG_0({\mbomega})},\\ \nonumber
\end{eqnarray}
and c is a constant. Let $P^{*}$ be the probability measure induced by $G^{*}$ on the
parameter space ${\mathcal N}$, and ${\bf u}, {\bf v} \in \partial {\mathcal A}$. Also let
${\mathcal L}=\{{\bf z}|~{\bf z}=\rho {\bf u}+(1-\rho){\bf v}, -\infty < \rho < \infty\}$.
If $P^{*}\{{\mbomega}^{'}{\bf(u-v)}\ne 0\}\ne 0$,  then by the inequality
\begin{eqnarray}
\mbox{exp}(\rho m_1+(1-\rho) m_2) \leq \rho~ e^{m_1}+(1-\rho)~e^{m_2}, ~0<\rho~<1,
\end{eqnarray}
for all real $m_1 $ and $m_2$ with equality if and only if $m_1=m_2$. Thus $t({\bf z}) < c$.
And hence the line ${\mathcal L}$ is not a support of ${\mathcal A}$ unless the intersection
of ${\mathcal L}$ and the closure of ${\mathcal A}$ is a single boundary point. If
$P^{*}\{{\mbomega}^{'}{\bf (u-v)} \ne 0\}=0$, then
\begin{eqnarray}
t({\bf z})=\int_{{\mathcal N}}\mbox{exp}\{{\mbomega}^{'}{\bf v}+\rho
{\mbomega}^{'}{\bf (u-v)}\}dG^{*}({\mbomega})=c.
\end{eqnarray}
Thus ${\mathcal L} \subseteq \partial {\mathcal A}$. ~~~~~~~~~~~~~Q.E.D.

\indent Obviously, Proposition 2.1 can be extended to the more general set-up, though it is
designed for the problem (1.4) in this paper. Let $P =\{1,\ldots, p\}$, and for every
$a$: {\reg}, let $a'$ be its complement and $|a|$ its cardinality. For each $a$, we partition
${\mXbar}_{n}$ and ${\bf S}_{n}$ as
\begin{eqnarray}
  {\mXbar}_{n}=\left(
  \begin{array}{c}
    {\mXbar}_{na}\\
    {\mXbar}_{na'}
  \end{array}
  \right) ~~\hbox{and}~~
  {\bf S}_{n} = \left(
  \begin{array}{cc}
    {\bf S}_{naa}& {\bf S}_{naa'}\\
    {\bf S}_{na'a}& {\bf S}_{na'a'}
  \end{array}
  \right),
\end{eqnarray}
and write
\begin{eqnarray}
  {\mXbar}_{na:a'} = {\mXbar}_{na}-{\bf S}_{naa'}{\bf S}_{na'a'}^{-1}
  {\mXbar}_{na'}, \\
  {\bf S}_{naa:a'}={\bf S}_{naa}-{\bf S}_{naa'}{\bf S}_{na'a'}^{-1}
  {\bf S}_{na'a}.
\end{eqnarray}
Further, let
\begin{eqnarray}
  I_{na} = 1\{{\mXbar}_{na:a'}>{\bf 0}, {\bf S}_{na'a'}^{-1}{\mXbar}_{na'}\leq 0\},
\end{eqnarray}
for {\reg}, where $1\{\cdot\}$ denotes the indicator function.
Then for the problem of testing $H_{0}: {\mbtheta} = {\bf 0}$ against the positive orthant
space $H_{1\it{O}^{+}}: {\mbtheta} \in {\cal O}^{+}_{p} \backslash \{{\bf 0}\}$, from the
results of Perlman (1969) and Sen and Tsai (1999) the LRT and the UIT statistics are of the
forms
\begin{eqnarray}
  L_{n} = \sum_{\mbox{\reg}}\left\{ {\frac{n{\mXbar}'_{na:a'}
  {\bf S}_{naa:a'}^{-1} {\mXbar}_{na:a'}} {1+n{\mXbar}'_{na'} {\bf S}_{na'a'}^
  {-1} {\mXbar}_{na'}}} \right\}I_{na}
\end{eqnarray}
and
\begin{eqnarray}
  U_{n}= \sum_{\mbox{\reg}}\{n{\mXbar}'_{na:a'}{\bf S}_{naa:a'}^
 {-1}{\mXbar}_{na:a'}\}I_{na}
\end{eqnarray}
respectively.

\indent Side by side, for testing $H_{0}: {\mbtheta} = {\bf 0}$ against the halfspace
$H_{1\it{H}^{*}}: {\mbtheta} \in {\cal H}^{*}_{p}\backslash \{{\bf 0}\}$ the LR and the UI
test statistics are of the forms
\begin{eqnarray}
  L_{n}^{\ast} = \sum_{\mbox{\rag}}\left\{ {\frac{n{\mXbar}'_{na:a'}
  {\bf S}_{naa:a'}^{-1} {\mXbar}_{na:a'}} {1+n{\mXbar}'_{na'} {\bf S}_{na'a'}^
  {-1} {\mXbar}_{na'}}} \right\}I^{*}_{na}
\end{eqnarray}
and
\begin{eqnarray}
  U_{n}^{\ast}= \sum_{\mbox{\rag}}\{n{\mXbar}'_{na:a'}{\bf S}_{naa:a'}^
 {-1}{\mXbar}_{na:a'}\}I^{*}_{na}
\end{eqnarray}
respectively, where $P_{1}=P-\{ p\}, I^{*}_{nP}=1\{\Bar{X}_{np}>0\}$ and
$I^{*}_{nP1}=1\{\Bar{X}_{np}\leq 0\}$.

\indent It is clear from (2.13)-(2.16) that
\begin{eqnarray}
 L_{n}^{\ast}\geq L_{n}, ~~\mbox{and}~~ U_{n}^{\ast}\geq U_{n},
\end{eqnarray}
with probability one.

\vspace{0.3cm}
\indent {\sc Theorem 2.1}. {\it For the problem of testing $H_{0}:{\mbtheta} = {\bf 0}$
against $H_{1}:{\mbtheta} \in {\cal C}\backslash \{{\bf 0}\}$, where ${\cal C}$ is either
the positive orthant space ${\cal O}_{p}^{+}$ or the halfspace ${\cal H}_{p}^{*}$, the
UIT is not a Bayes test.}

\indent {\sc Proof}. The UIT statistics in (2.14) and in (2.16) have a similar
structure, hence without loss of generality it suffices to consider the problem of testing
against the positive orthant space. Let $u_{\alpha}$ be the critical point of the level of
significance $\alpha$ for the UIT, which can be obtained from the right hand side of (3.9).
Let ${\mathcal A}_{U}=\cup_{\emptyset}^{P}{\cal A}_{Ua}$ be the acceptance region of the UIT
for testing $H_{0}: {\mbtheta} = {\bf 0}$ against
$H_{1}: {\mbtheta} \in {\cal O}^{+}_{p} \backslash \{{\bf 0}\}$,
where ${\cal A}_{Ua}=\{({\mXbar}_{n}, {\bf S}_{n})| ~n{\mXbar}'_{na:a'}
{\bf S}^{-1}_{naa:a'}{\mXbar}_{na:a'} \leq u_{\alpha}\}I_{na},$ {\reg}.
For the situation $|a|=1$, we may note that there exists a support ${\mathcal L}_{a}$ of the
set ${\cal A}_{U}$ such that $({\mathcal L}_{a}\cap I_{na})\subseteq \partial{\cal A}_{Ua}$
and ${\mathcal L}_{a} \nsubseteq \partial{\cal A}_{U}$, since ${\mathcal L}_{a} \cap I_{na}$
is easily seen to be a half-line. Thus theorem follows by Proposition 2.1. ~~~~~~~~~~~~~Q.E.D.

\indent Note that for the problem of testing against the positive orthant space or the problem
of testing aginst the half-space, the LRT is isomorphic to the UIT when the covariance matrix
is known $({\mbSigma}={\mbSigma}_{0}$, with ${\bf S}_{n}$ being replaced by
${\mbSigma}_{0}$) or is partially unknown, (${\mbSigma}={\sigma}^{2}{\mbSigma}_{0}$, where
${\sigma}^{2}$ is an unknown scalar and ${\mbSigma}_{0}$ is a known positive definite matrix,
with ${\bf S}_{n}$ being replaced by $s^{2}{\mbSigma}_{0}$). Hence, similar arguments as in
the proof of theorem 2.1, we may conclude that the LRT is not a Bayes when the covariance
matrix is known or partially unknown for the problems of testing against restricted
alternatives. As such, when the covraiance matrix is totally unknown the LRT is asymptotic
equivalent to the UIT as the sample size is sufficiently large. Thus, the LRT is also not a
Bayes test when sample size is sufficiently large for the problems of testing against
restricted alternatives..

\indent Next, consider for the finite union-intersection test (FUIT, Roy et al., 1972) based on the one-sided
coordinatedwise Student t-tests. Define ${\bf S}_{n}=(S_{nij})$ as in (2.9) and let
\begin{eqnarray}
t_j=\frac{{\sqrt {n}}\bar X_{nj}}{\sqrt {S_{njj}}},~ j=1, \cdots, p.
\end{eqnarray}
Corresponding to a given significance level $\alpha$, define
\begin{eqnarray}
{\alpha}^{*}:~p{\alpha}^{*}=\alpha.
\end{eqnarray}
Let $t_{n-1,{\alpha}^{*}}$ be the upper
$100{\alpha}^{*}$\% point of the Studentn t-distribution with $n-1$ degrees of freedom.
Consider the critical region ${\cal W}_j=\{t_j|~t_j\geq t_{n-1,{\alpha}^{*}}\}$ for
$~j=1,\cdots, p$. Then the critical region of the FUIT is
\begin{eqnarray}
{\cal W}^{*}=\cup_{j=1}^{p}{\cal }{\cal W}_j
\end{eqnarray}
and the acceptance region is
\begin{eqnarray}
{\cal A}^{*}=\cap_{j=1}^{p}{\cal }\Bar{\cal A}_j,
\end{eqnarray}
where ${\cal A}_j=R\setminus{\cal W}_j$ and $\Bar{\cal A}_j$ denote the closure of
${\cal A}_j, j=1, \cdots, p$. Therefore, by Proposition 2.1 we have the following.

\vspace{0.3cm}
\indent {\sc Theorem 2.2}. {\it For the problem of testing $H_{0}:{\mbtheta} = {\bf 0}$
against $H_{1}:{\mbtheta} \in {\cal C}\backslash \{{\bf 0}\}$, where ${\cal C}$ is either
the positive orthant space ${\cal O}_{p}^{+}$ or the halfspace ${\cal H}_{p}^{*}$, the
FUIT is not a Bayes test.}

\indent In passing, we note that the explicit forms of the Bayes tests are hardly obtained
due to the difficulty of integration over restricted parameter spaces. Hence, for the
problems of testing against restricted alternatives (1.4) we will study the Bayesian tests
only in terms of the non-informative prior distribution for the nuisance parameter
${\mbSigma}$
\begin{eqnarray}
   g({\mbSigma}) \propto  ~|{\mbSigma}|^{-\frac {1}{ 2}(p+1)} ~d\nu({\mbSigma}),
            ~ {\mbSigma}\in {\mathcal M}
\end{eqnarray}
for some measure $\nu(\cdot)$, where ${\mathcal M}$ is the set of positive definite matrices.
First, consider the unique Haar measure in ${\mbSigma}$ over the set
${\cal M}$, i.e., $\nu({\mbSigma})={\mbSigma}$. Then the marginal density of ${\mXbar}_{n}$
and ${\bf S}_{n}$ is given by
\begin{eqnarray}
  P_{1}({\mXbar}_{n}, {\bf S}_{n}|~{\mbtheta}) \propto ~|{\bf S}_{n}|^{\frac{1}{2}(n-p-2)}
  |{\bf S}_{n}+n({\mXbar}_{n}-{\mbtheta})({\mXbar}_{n}-{\mbtheta})'|^{-\frac{1}{2}(n-1)},
\end{eqnarray}
namely,
\begin{eqnarray}
  P_{1}({\mXbar}_{n}, {\bf S}_{n}|~{\mbtheta}) \propto ~|{\bf S}_{n}|^{-\frac{1}{2}(p+1)}
   [1+n({\mXbar}_{n}-{\mbtheta})'{\bf S}_{n}^{-1}
   ({\mXbar}_{n}-{\mbtheta})]^{-\frac{1}{2}(n-1)}.
\end{eqnarray}
Based on the density in (2.20), then
\begin{eqnarray}
&~&\frac{\sup_{\mbtheta \in {\cal O}_{p}^{+}}P_{1}({\mbtheta} | {\mXbar}_{n}, {\bf S}_{n})}
{P_{1}({\bf 0} | {\mXbar}_{n}, {\bf S}_{n})} \\
  & = & \left[ \frac{1+n{\mXbar}_{n}'{\bf S}_{n}^{-1}{\mXbar}_{n}}
    {1+\inf_{\mbtheta \in {\cal O}^{+}_{p}}n({\mXbar}_{n}-{\mbtheta})'{\bf S}_{n}^{-1}
     ({\mXbar}_{n}-{\mbtheta})} \right ]^{\frac{1}{2}(n-1)} \nonumber \\
 & = & \left [1+|\!| {\pi}_{{\bf S}_{n}}(n^{1/2} {\mXbar}_{n}; {\cal O}^{+}_{p})
    |\!|_{{\bf S}_{n}}^{2} \{ 1 + |\!| n^{1/2} {\mXbar}_{n} - {\pi}_{{\bf S}_{n}}
     (n^{1/2} {\mXbar}_{n}; {\cal O}^{+}_{p} ) |\!|_{{\bf S}_{n}}^{2} \}^{-1}\right
     ]^{\frac{1}{2}(n-1)} \nonumber \\
  & = & (1+ L_{n})^{\frac{1}{2}(n-1)}. \nonumber
\end{eqnarray}

\indent For any ${\bf b} \in R^{p}\backslash\{{\bf 0}\}$, the following density can be
obtained from (2.19)
\begin{eqnarray}
  P_{2}({\mXbar}_{n}, {\bf S}_{n}, {\bf b}|~{\mbtheta}) &\propto& ~|{\bf b}'{\bf S}_{n}
  {\bf b}|^{-\frac{1}{2}(p+1)} [1+n {\bf b}'({\mXbar}_{n}-{\mbtheta})
  ({\mXbar}_{n}-{\mbtheta})'{\bf b}/{\bf b}'{\bf S}_{n}{\bf b}]^{-\frac{1}{2}(n-1)}.
\end{eqnarray}
Based on the density in (2.22), for each ${\bf b} \in R^{p}\backslash\{{\bf 0}\}$, we have
\begin{eqnarray}
&~&\frac{\sup_{\mbtheta \in {\cal O}_{p}^{+}}P_{2}({\mbtheta} | {\mXbar}_{n},
{\bf S}_{n}, {\bf b})}
{P_{2}({\bf 0} | {\mXbar}_{n}, {\bf S}_{n}, {\bf b})} \\
  & = & \left[ \frac{1+n {\bf b}'{\mXbar}_{n}{\mXbar}'_{n}{\bf b}/
    {\bf b}'{\bf S}_{n} {\bf b}}{1 + \inf_{\mbtheta \in {\cal O}_{p}^{+}}
    n({\mXbar}_{n}-{\mbtheta})'{\bf b} {\bf b}'({\mXbar}_{n}-{\mbtheta})/
    {\bf b}' {\bf S}_{n}{\bf b}} \right ]^{\frac{1}{2}(n-1)}. \nonumber
\end{eqnarray}
Partition ${\bf b}'=({\bf b}'_{a}, {\bf b}'_{a'})$ the same as in (2.9) and for each
$a$, \reg, let ${\bf G}_{a}= \left[
\begin{array}{cc}
    {\bf I}_{a}& {\bf -S}_{naa'}{\bf S}^{-1}_{na'a'}\\
    {\bf 0}& {\bf I}_{a'}
\end{array}
  \right].$
Further multiply ${\bf G}_{a}$ on the left of ${\mbtheta}$ and ${\mXbar}_{n}$
and its inverse on the right of ${\bf b}^{'}$. Next by using Lemma 1.1 of N\"{u}esch
(1966) and choosing ${\bf b}_{a'}= - {\bf S}_{na'a'}{\bf S}_{na'a}{\bf b}_{a}$,
then after some simplifications the $ r.h.s.$ of (2.27) reduces to
\begin{eqnarray}
\left[ 1+ \sum_{\mbox{\reg}} \sup_{{\bf b}_{a} \neq {\bf 0}}\left\{ {\frac{n({\bf b}'_{a}
{\mXbar}_{na:a'})^{2}}
{{\bf b}'_{a}{\bf S}_{naa:a'}{\bf b}_{a}}} \right\}I_{na} \right]^{\frac{1}{2}(n-1)}.
\end{eqnarray}
Since (2.28) holds for all ${\bf b}_{a} \neq {\bf 0}$, \reg, further note that
\begin{eqnarray}
&~&\sum_{\mbox{\reg}}\sup_{{\bf b}_{a} \neq {\bf 0}}\left\{ {\frac{n({\bf b}'_{a}
        {\mXbar}_{na:a'})^{2}}{{\bf b}'_{a}{\bf S}_{naa:a'}{\bf b}_{a}}} \right\}I_{na}\\
& = &  \sum_{\mbox{\reg}}\{n{\mXbar}'_{na:a'}{\bf S}_{naa:a'}^
         {-1}{\mXbar}_{na:a'}\}I_{na} \nonumber \\
& = & U_{n}. \nonumber
\end{eqnarray}

\indent If a Bernardo reference prior (Chang and Eaves, 1990) of the form
\begin{eqnarray}
  g_{1}({\mbSigma}) \propto ~ |{\mbSigma}|^{-\frac{1}{2}(p+1)}
     |{\bf I}+{\mbSigma}*{\mbSigma}^{-1}|^{-\frac{1}{2}}~d{\mbSigma},
     ~{\mbSigma}\in {\mathcal M}
\end{eqnarray}
is considered, where ${\bf C}*{\bf D}$ denotes the Hadamard product of the matrices
${\bf C}$ and ${\bf D}$, then the marginal density is
\begin{eqnarray}
  h({\bf W}_{n})= c \int_{{\mbSigma}\in {\mathcal M}}
               |{\mbSigma}|^{-\frac{1}{2}(p+1)}|{\bf I}+
         {\mbSigma}*{\mbSigma}^{-1}|^{-\frac{1}{2}}~
          \mbox{e}^{{-\frac{1}{2}}\mbox{tr}({\mbSigma}^{-1}
          {\bf W}_{n})}~d{\mbSigma},
\end{eqnarray}
where ${\bf W}_{n}={\bf S}_{n}+n({\mXbar}_{n}-{\mbtheta})({\mXbar}_{n}-{\mbtheta})'$ and
$c$ is the norming constant. By the facts that
$d|{\bf W}_{n}|= |{\bf W}_{n}|\mbox{tr}({\bf W}^{-1}_{n})(d{\bf W}_{n})$ and
$d\mbox{tr}({\mbSigma}^{-1}{\bf W}_{n})=\mbox{tr}({\mbSigma}^{-1})(d{\bf W}_{n})$, where
$(d{\bf W}_{n})$ denotes the exterior differential form of ${\bf W}_{n}$. Then
\begin{eqnarray}
\frac{dh({\bf W}_{n})}{d|{\bf W}_{n}|}
    &=& {\frac{-c}{2}}\int_{{\mbSigma}\in {\mathcal M}}|{\bf W}_{n}|^{-1}
         (\mbox{tr}{\bf W}^{-1}_{n})^{-1}\mbox{tr}({\mbSigma}^{-1})
         |{\mbSigma}|^{-\frac{1}{2}(p+1)}  \\
     & & \hspace{3cm} \times |{\bf I}+{\mbSigma}*{\mbSigma}^{-1}|^{-\frac{1}{2}}
         \mbox{e}^{{-\frac{1}{2}}\mbox{tr}({\mbSigma}^{-1}{\bf W}_{n})}
         ~d{\mbSigma} \nonumber  \\
     &< &  0 , \nonumber
\end{eqnarray}
hence $h({\bf W}_{n})$ is strictly decreasing in $|{\bf W}_{n}|$. Thus the corresponding
integrated LRT and integrated UIT can be constructed based on $|{\bf W}_{n}|$ only.

\indent Next, we consider the case if ${\mbSigma}$ is assigned a prior distrbution
$W^{-1}_{p}(\mbGamma, m)$, where $W^{-1}_{p}(\mbGamma, m)$ denoting the $p$-dimensional
inverted Wishart distribution with $m$ degrees of freedom and expectation $m{\mbGamma}$
which is positive definite. Then by Theorem 7.7.2 of Anderson (2003), the marginal density
is
\begin{eqnarray}
  h_{1}({\mXbar}_{n}, {\bf S}_{n}|~{\mbtheta}) & \propto & ~|{\mbGamma}|^{\frac{m}{2}}
      |{\bf S}_{n}|^{\frac{1}{2}(n-p-2)}|{\bf S}_{n}+n({\mXbar}_{n}
      -{\mbtheta})({\mXbar}_{n} -{\mbtheta})^{'}+{\mbGamma}|^{-\frac{1}{2}(n+m-1)}   \\
    & = &  |{\mbGamma}|^{\frac{m}{2}}|{\bf S}_{n}|^{\frac{1}{2}(n-p-2)}
         |{\bf W}_{n}+{\mbGamma}|^{-\frac{1}{2}(n+m-1)}, \nonumber
\end{eqnarray}
which is a monotone function of $|{\bf W}_{n}|$.

\indent Similarly, by (2.25) and (2.29) the $L_{n}^{*}$ and $U_{n}^{*}$ can be exactly obtained
with respect to those non-informative prior differentials mentioned above over the space
${\mbSigma}\in {\mathcal M}$ for the problem of testing against a halfspace ${\cal H}_{p}^{*}$.

\indent Let $\mathcal P_{p}$ be the space of positive definite matrices. Riemannian geometry
yields an invariant volume element $dv$ on $\mathcal P_{p}$ which is of the form
$dv=(\mbox{det}{\mbSigma})^{-(p+1)/2}(d{\mbSigma})$, where
$(d{\mbSigma})=\prod_{1 \leq j \leq i  \leq p} d{\sigma}_{ij}$ with $d{\sigma}_{ij}$ being
the Lebesgue measure on $R$. We may naturally adopt this relevant Haar measure $dv$ as the
underline measure of the a prior of $\mbSigma$. The likelihood integrated with respect to this
relevant Haar measure $dv$ for the nuisance covariace matrix can be viewed as in a Bayesian
light. It is not hard to note that for testing against global alternative, the Hotelling's $T^{2}$-test
is the version of integrated LRT. Moreover, for the problem (1.4) the corresponding integrated
LRT and integrated UIT also have the same forms as those of the corresponding LRT and UIT,
respectively. These will reconcile frequentist with Bayesian paradigms via
the fiducial inference for the problem (1.4).

\vspace{0.3cm}
\indent {\sc Theorem 2.3}. {\it For the problem of testing $H_{0}:{\mbtheta} = {\bf 0}$
against $H_{1}:{\mbtheta} \in {\cal C}\backslash \{{\bf 0}\}$, where ${\cal C}$ is either
the positive orthant space ${\cal O}_{p}^{+}$ or the halfspace
${\cal H}_{p}^{*}=\{{\mbtheta}|~{\theta}_{p} \geq 0 \}$, then both the LRT and the UIT
are the versions of integrated LRT and integrated UIT respectively.}

\vspace{0.3cm}
\def \theequation{3.\arabic{equation}}
\setcounter{equation}{0}

\indent {\bf 3. Power comparision or reporting $P$-value?} Historically, when it comes to hypothesis
testing problems, two different approaches are adopted to evaluate the tests for frequientist
paradigm: one is Fisher approcah to report $p$-value of the test, the other is Neyman-Pearson
approach to carry out the power of test at a fixed significance level. Both Fisher approach
and Neyman-Pearson approach are well adopted to measure the performance of the Hotelling's
$T^{2}$-test for the problem of testing against the unrestricted alternatives. We would like
to know whether the Fisher approach and/or the Neyman-Pearson approach can be well adopted for
the problem (1.4). First, consider the problem of testing against the half-space.

\indent Let $\chi_{m}^{2}$ be the central chisquare distribution with $m(\geq 0)$ degrees of
freedom. Denoted by
\begin{eqnarray}
  G_{a,b}(u)=P\{\chi_{a}^{2}/\chi_{b}^{2}\leq u\}, ~u\in R^{+},
\end{eqnarray}
and ${\bar G}_{a,b}(u)=1-G_{a,b}(u)$, where $R^{+}$ denotes the positive real number.  Also
let
\begin{eqnarray}
  {\bar G}_{n,a,p}^{\ast}(u) = \int_{0}^{\infty} {\bar G}_{a,n-p}({\frac{u}{1+t
  }})dG_{p-a,n-p+a}(t),~~u\in R^{+},
\end{eqnarray}
which is the convolution of the d.f.'s of $\chi_{a}^{2}/\chi_{n-p}^{
2}$ and 1+$\chi_{p-a}^{2}/\chi_{n-p+a}^{2}$. By (2.15) and (2.16), as such for the problem
of testing against the halfspace ${\mathcal H}^{*}_{p} \backslash \{{\bf 0}\}$, it is easy
to show that when ${\mbSigma} \in {\cal M},~\forall ~c > 0$
\begin{eqnarray}
 P_{{\bf 0}, {\mbSigma}}\{L_{n}^{\ast}\geq c\}
 &=& {\frac{1}{2}}\left [\bar G_{p-1,n-p}(c)+\bar G_{p,n-p}(c)\right],
\end{eqnarray}
and
\begin{eqnarray}
 P_{{\bf 0}, {\mbSigma}}\{U_{n}^{\ast}\geq c\}
  &=&{\frac{1}{2}}\left [\bar G_{n,p-1,p}^{*}(c)+\bar G_{n,p,p}^{*}(c)\right].
\end{eqnarray}
Notice that (3.3) and (3.4) are free from the nuisance parameter ${\mbSigma}$, namely,
the $p$-values of LRT and UIT do not depend on the nuisance parameter. By the result of
Tang (1994) and similar arguments in Theorems 2.1 and 3.1 of Sen and Tsai (1999) we then
obtain that both LRT and UIT for testing against a halfspace are similar and unbiased.
Hence, for the problem of testing against a halfspace space both approaches of Fisher and
Neyman-Pearson can well explain the performance of the LRT and the UIT.

\indent As for the problem of testing against the positive orthant space, it follows from
the arguments of Perlman (1969) and Sen and Tsai (1999) that for every c $>$ 0,
\begin{eqnarray}
  P\{ L_{n} \geq c| H_0,{\mbSigma} \} = P_{{\bf 0},{\mbSigma}} \{ L_{n} \geq c \} \\
  = \sum_{k=1}^{p}w(p,k;{\mbSigma})P\{ \chi_{k}^{2} /\chi_{n-p}^{2}\geq c\}\nonumber
\end{eqnarray}
and
\begin{eqnarray}
  P\{ U_{n} \geq c| H_0,{\mbSigma} \} = P_{{\bf 0},{\mbSigma}} \{ U_{n} \geq c \} \\
  = \sum_{k=1}^{p}w(p,k;{\mbSigma}) {\bar G}_{n,k,p}^{\ast}(c),\nonumber
\end{eqnarray}
respectively, where for each $k$ (=0,1,\ldots,$p$)
\begin{eqnarray}
  w(p,k;{\mbSigma})=\sum\nolimits_{\{a\subseteq P: |a|=k\}}P\{{\bf Z}_{a:a'}>{\bf 0}, {
  \mbSigma}_{a'a'}^{-1}{\bf Z}_{a'} \leq {\bf 0}\},
\end{eqnarray}
and ${\bf Z}\sim{\it N}_{p}({\bf 0},{\mbSigma})$ and the partitions are made as in
(2.9)-(2.11) (with $({\mXbar}_{n}, {\bf S}_{n})$ being replaced by ({\bf Z}, {\bSigma})).

\indent By virtue of (3.5), (3.6) and (3.7), we encounter the problem for which both the null
distributions of LRT and UIT statistics depend on the unknown covariance matrix, though
the dependence of (3.5) and (3.6) on {\bSigma} is only through the $w(p,k;{\mbSigma})$.
Hence, we have difficulty to adopt Fisher approach to report the $p$-value. On the other hand,
according to the definition of level of significance, Perlman (1969) and Sen and Tsai (1999)
allowed ${\mbSigma}$ to vary over the entire class ${\mathcal M}$ and obtained that for every
c $> 0$,
\begin{eqnarray}
  \sup_{\{{\mbSigma}\in {\mathcal M}\}}P_{{\bf 0},{\mbSigma}}\{L_{n}\geq c\}
  = {\frac {1}{2}}\left[\bar G_{p-1,n-p}(c)+\bar G_{p,n-p}(c)\right]
\end{eqnarray}
and
\begin{eqnarray}
  \mathop{\sup}_{\{{\mbSigma}\in {\mathcal M}\}}P_{{\bf 0},{\mbSigma}}\{U_{n} \geq
  c\}={\frac{1}{2}}\left [\bar G_{n,p-1,p}^{*}(c)+\bar G_{n,p,p}^{*}(c)\right],
\end{eqnarray}
respectively. If the right hand side of (3.8) and (3.9) are equated to $\alpha$, we then
get the critical values of LRT and UIT respectively. However, based on them we will
obtain the results which are against our common statistical sense. To proceed them,
some notations are needed.

\indent Let $\pi_{\bf A}({\bf x};{\cal C})$ be the orthogonal projection of ${\bf x}$ onto
${\cal C}$ with respect to the inner product $<, >_{\bf A}$, then
\begin{eqnarray}
   U_{n}=|\!|{\pi}_{{\bf S}_n}(n^{1/2}{\mXbar}_{n};{\cal O}_{p}^{+})
        |\!|_{{\bf S}_n}^{2}
\end{eqnarray}
and
\begin{eqnarray}
   U^{*}_{n}=|\!|{\pi}_{{\bf S}_n}(n^{1/2}{\mXbar}_{n}; {\mathcal H}^{*}_{p}\})
        |\!|_{{\bf S}_n}^{2}.
\end{eqnarray}

\vspace{0.3cm}
\indent {\sc Proposition 3.1.} {\it The UIT for testing $H_{0}:{\mbtheta}={\bf 0}$ against
$H_{1}: {\mbtheta} \in {\cal C}\backslash \{{\bf 0}\}$, where ${\cal C}$ is a closed convex cone
such that ${\cal U} \subseteq {\cal C} \subseteq {\mathcal H}^{*}_{p}$ with ${\cal U}$ being
a p-dimensional open set and being contained in the halfspace ${\mathcal H}^{*}_{p}$, is dominated by
the UIT for testing $H_{0}:{\mbtheta} ={\bf 0}~{against}~H_{1}: {\mbtheta}
\in {\mathcal H}^{*}_{p} \backslash \{{\bf 0}\}$.}

\indent {\sc Proof.} First note that ${\mathcal C} \subseteq {\mathcal H}^{*}_{p}$, hence
\begin{eqnarray}
 |\!| \pi_{{\bf S}_n}(n^{1/2}{\mXbar}_{n}; {\mathcal H}^{*}_{p}) |\!|_{{\bf S}_n}^{2}
  \geq |\!| \pi_{{\bf S}_n}(n^{1/2}{\mXbar}_{n}; {\cal C}) |\!|_{{\bf S}_n}^{2}. \nonumber
\end{eqnarray}
Then for every $c >0,$
\begin{eqnarray}
   P_{{\mbtheta},{\mbSigma}} \{|\!|{\pi}_{{\bf S}_n}(n^{1/2}{\mXbar}_{n};
   {\mathcal H}^{*}_{p})|\!|_{{\bf S}_n}^{2}  \geq c \} \geq P_{{\mbtheta},{\mbSigma}}
   \{ |\!| \pi_{{\bf S}_{n}}(n^{1/2}{\mXbar}_{n}; {\mathcal C}) |\!|_{{\bf S}_{n}}^{2}
   \geq c \}. \nonumber
\end{eqnarray}
\indent Let ${\mbSigma}^{-1}
={\bf B}^{'}{\bf B}$, and ${\mYbar}_{n}={\bf B}{\mXbar}_{n}$, and ${\bf T}_{n}
={\bf B}{\bf S}_{n}{\bf B}^{'}.$ Since ${\bf B}{\mathcal H}^{*}_{p}={\mathcal H}^{+}_{p}$
(another halfspace), thus $P_{{\bf 0},{\mbSigma}} \{|\!|{\pi}_{{\bf S}_n}(n^{1/2}
{\mXbar}_{n}; {\mathcal H}_{p}^{+})|\!|_{{\bf S}_n}^{2} \geq c \} = P_{{\bf 0},{\bf I}}
\{|\!|{\pi}_{{\bf T}_n}(n^{1/2}{\mYbar}_{n};{\mathcal H}_{p}^{+})|\!|_{{\bf T}_n
  }^{2} \geq c \}$. Notice that ${\cal E}({\bf Y}_{n})=$(n-1)${\bf I}$,
thus $P_{{\bf 0},{\mbSigma}} \{|\!|{\pi}_{{\bf S}_n}(n^{1/2} {\mXbar}_{n};
{\mathcal H}_{p}^{+})|\!|_{{\bf S}_n}^{2} \geq c \}$ is independent of the unknown
${\mbSigma}$. Consider the sequence $\{{\mbSigma}_{k}\}$ and let ${\mbSigma}_{k}^{-1}
={\bf B}_{k}^{'}{\bf B}_{k}$, for all $k \geq 1$, such that ${\bf B}_{k}
({\mathcal C_{\lambda}}) \subseteq {\bf B}_{k+1}({\cal C_{\lambda}})$ and
$\cup_{k=1}^{\infty}{\bf B}_{k} ({\mathcal C_{\lambda}}) =int({\cal H}_{p}^{+})$, where
${\mathcal C_{\lambda}}$ with $0 <{\lambda} < 1$ is a right circular cone defined as in
$2.3^{o}$ of Perlman (1969). The same proof as in Theorem 6.2 of Perlman (1969), we obtain
\begin{eqnarray}
 & &P_{{\bf 0},{\mbSigma}_{k}} \{|\!|{\pi}_{{\bf S}_n}(n^{1/2}{\mXbar}_{n};
{\mathcal H}_{p}^{+})|\!|_{{\bf S}_n }^{2} \geq c \} \geq
 P_{{\bf 0},{\mbSigma}_{k}}\{ |\!| \pi_{{\bf S}_{n}} (n^{1/2}{\mXbar}_{n};
  {\mathcal C}) |\!|_{{\bf S}_{n}}^{2} \geq c \}\nonumber \\
   & &\geq  P_{{\bf 0},{\mbSigma}_{k}}\{ |\!| \pi_{{\bf S}_{n}} (n^{1/2}{\mXbar}_{n};
  {\mathcal C_{\lambda}}) |\!|_{{\bf S}_{n}}^{2} \geq c \} =
  P_{{\bf 0},{\bf I}}\{ |\!| \pi_{{\bf S}_{n}^{*}} (n^{1/2}{\mXbar}_{n}^{*};
  {\bf B}_{k}({\mathcal C_{\lambda}})) |\!|_{{\bf S}_{n}^{*}}^{2} \geq c \},\nonumber
\end{eqnarray}
where ${\mXbar}_{n}^{*}={\bf B}_{k}{\mXbar}_{n}$, and ${\bf S}^{*}_{n}=
{\bf B}_{k}{\bf S}_{n}{\bf B}_{k}^{'}$. By Lemma 2.8 of Brown (1986), we may conclude that
the power functions of the UIT for testing $H_{0}:{\mbtheta}={\bf 0}~\hbox{vs.}~H_{1}:
{\mbtheta} \in {\mathcal C}$ are analytic on the space $\{{\mbtheta} \in {\mathcal C},
{\mbSigma} > {\bf 0}\}$, and hence they are continuous. Thus
$\lim_{k \to \infty}P_{{\bf 0},{\bf I}}\{ |\!| \pi_{{\bf S}^{*}_{n}}
   (n^{1/2}{\mXbar}_{n}^{*}; {\bf B}_{k}({\cal C_{\lambda}})) |\!|_{{\bf S}^{*}_{n}}^{2}
   \geq c \}= P_{{\bf 0},{\bf I}}\{ |\!| \pi_{{\bf S}_{n}^{*}} (n^{1/2}{\mXbar}_{n}^{*};
    {\cal H}_{p}^{+}) |\!|_{{\bf S}^{*}_{n}}^{2} \geq c \}$, and hence
\begin{eqnarray}
{\sup}_{\{{\mbSigma} \in
{\cal M}\}} P_{{\bf 0},{\mbSigma}}\{ |\!| \pi_{{\bf S}_{n}} (n^{1/2}{\mXbar}_{n};
  {\cal C}) |\!|_{{\bf S}_{n}}^{2} \geq c \}
  =P_{{\bf 0},{\mbSigma}} \{|\!|{\pi}_{{\bf S}_n}(n^{1/2}{\mXbar}_{n};
    {\cal H}_{p}^{+})|\!|_{{\bf S}_n }^{2} \geq c \}.
\end{eqnarray}
Thus, the proposition follows. ~~~~~~~~~~~~~Q.E.D.

\indent As for the LRT, first, note that the LRT statistics in (2.13) and (2.15) can be
represented as
\begin{eqnarray}
  L_{n}= |\!| {\pi}_{{\bf S}_{n}}(n^{1/2} {\mXbar}_{n}; {\cal O}^{+}_{p})
    |\!|_{{\bf S}_{n}}^{2} \{ 1 + |\!| n^{1/2} {\mXbar}_{n} - {\pi}_{{\bf S}_{n}}
     (n^{1/2} {\mXbar}_{n}; {\cal O}^{+}_{p} ) |\!|_{{\bf S}_{n}}^{2} \}^{-1}
\end{eqnarray}
and
\begin{eqnarray}
   L^{*}_{n}=|\!| {\pi}_{{\bf S}_{n}}(n^{1/2} {\mXbar}_{n}; {\cal H}^{*}_{p})
    |\!|_{{\bf S}_{n}}^{2} \{ 1 + |\!| n^{1/2} {\mXbar}_{n} - {\pi}_{{\bf S}_{n}}
     (n^{1/2} {\mXbar}_{n}; {\cal H}^{*}_{p} ) |\!|_{{\bf S}_{n}}^{2} \}^{-1}
\end{eqnarray}
respectively. Also note that for any ${\mathcal C} \subseteq {\mathcal H}^{*}_{p}$, $
 |\!| \pi_{{\bf S}_n}(n^{1/2}{\mXbar}_{n}; {\mathcal H}^{*}_{p}) |\!|_{{\bf S}_n}^{2}
  \geq |\!| \pi_{{\bf S}_n}(n^{1/2}{\mXbar}_{n}; {\cal C}) |\!|_{{\bf S}_n}^{2}$, and hence
$|\!| n^{1/2} {\mXbar}_{n} - {\pi}_{{\bf S}_{n}}(n^{1/2} {\mXbar}_{n}; {\cal H}^{*}_{p} )
|\!|_{{\bf S}_{n}}^{2} \leq |\!| n^{1/2} {\mXbar}_{n} - {\pi}_{{\bf S}_{n}}
     (n^{1/2} {\mXbar}_{n}; {\cal C}) |\!|_{{\bf S}_{n}}^{2}$. Thus, by similar argument
as in the proof of proposition 3.1, we also have the following.

\vspace{0.3cm}
\indent {\sc Proposition 3.2.} {\it The LRT for testing $H_{0}:{\mbtheta}={\bf 0}~against ~H_{1}: {\mbtheta}
\in {\cal C}\backslash \{{\bf 0}\}$, where ${\cal C}$ is a closed convex cone such that
${\cal U} \subseteq {\cal C} \subseteq {\mathcal H}^{*}_{p}$ with ${\cal U}$ being a p-dimensional open set
and being contained in the halfspace ${\mathcal H}^{*}_{p}$, is dominated by the LRT for testing
$H_{0}:{\mbtheta} ={\bf 0}~{against}~H_{1}: {\mbtheta} \in {\mathcal H}^{*}_{p} \backslash \{{\bf 0}\}$.}

\indent Proposition 3.2 implies the phenomenon first noted by Tang (1994) that the LRT for
testing against the positive orthant space is dominated by that of testing against the
halfspace. Note the propositions 3.1 and 3.2 are no longer to be true when the sample size is
sufficiently large.

\indent To overcome these unpleasant phenomena, instead of using the conservative maximum
principle to define the level of significance $\alpha$, i.e., taking the supremum of rejected
probability under null hypothesis over the set ${\cal M}$, we may use the Bayesian notions, by
taking the average of rejected probability under null hypothesis over the set ${\cal M}$ with
respect to the weight function $g({\mbSigma})$ of ${\mbSigma}$, to define the level of
significance $\alpha$.

\indent The joint density of $({\mXbar}_{n},{\bf S}_{n})$ under the null hypothesis is
\begin{eqnarray}
f({\mXbar}_{n}, {\bf S}_{n})=k_{0}|{\mbSigma}|^{-\frac{1}{2}n}
         |{\bf S}_{n}|^{\frac{1}{2}(n-p-2)}
         \mbox{e}^{{-\frac{1}{2}}\mbox{tr}({\mbSigma}^{-1}{\bf V}_{n})},
\end{eqnarray}
where $k_{0}$ is the norming constant and
${\bf V}_{n}={\bf S}_{n}+n{\mXbar}_{n}{\mXbar}^{'}_{n}$.
First consider the inverted Wishart distribution $W^{-1}({\mbGamma},m)$ (Anderson, 2003),
which is a proper prior of ${\mbSigma}$, as the weight function. Then we have
\begin{eqnarray}
 h_{1}({\mXbar}_{n},{\bf S}_{n})&=&\int_{{\mbSigma}\in {\mathcal M}}
     f({\mXbar}_{n}, {\bf S}_{n})g({\mbSigma}) ~d{\mbSigma} \\ \nonumber
 &=&  k_{2}~|{\mbGamma}|^{\frac{1}{2}m}|{\bf S}_{n}|^{\frac{1}{2}(n-p-2)}
        |{\bf V}_{n}+{\mbGamma}|^{-\frac{1}{2}(n+m-1)},
\end{eqnarray}
where the constant $k_{2}$ can be determinated by the equation
$\int_{{\cal D}} ~h_{1}({\mXbar}_{n},{\bf S}_{n})~~dx_{1}\cdots dx_{p}=1$ with
${\cal D}=\{(x_{1},\cdots,x_{p})|~{\bf S}_{n}\in {\mathcal M} ~\hbox{in probability}\}$.
Let
\begin{eqnarray}
  b_{1}(k,n,p)=\sum_{|a|=k}\int_{I_{na}\cap {\cal D}}~h_{1}({\mXbar}_{n},{\bf S}_{n})
   ~dx_{1}\cdots dx_{p},
\end{eqnarray}
then it is easy to note that $ b_{1}(k,n,p) > 0, \forall k=1, 2, \cdots,p$ and
$\sum_{k=0}^{p}b_{1}(k,n,p)=1$.
Thus by Fubini theorem the critical values of the  LRT and the UIT for testing against the
positive orthant space can then be determinated via the formula
\begin{eqnarray}
  \alpha=\sum_{k=0}^{p}b_{1}(k,n,p){\bar G}_{k,n-p}(d^{1}_{n,{\alpha}})
\end{eqnarray}
and
\begin{eqnarray}
  \alpha=\sum_{k=0}^{p}b_{1}(k,n,p){\bar G}_{n,k,p}^{\ast}(c^{1}_{n,{\alpha}})
\end{eqnarray}
respectively. If the critical values of LRT and UIT are determined by the equations (3.18) and
(3.19), then obviously both the LRT and the UIT are similiar and unbiased. These conclusions
are contrary to the ones made by Perlman (1969) and Sen and Tsai (1999) that both the LRT and
the UIT are not similar and biased.

\indent The level of significance $\alpha$ may be defined by using the noninformative priors
such as the Haar measure or the Bernardo reference prior, in which the corresponding posterior
densities exist, as the weight functions. Suppose that the Haar measure of ${\mbSigma}$ is
taken as the weight function, then we have
\begin{eqnarray}
  h_{2}({\mXbar}_{n},{\bf S}_{n})&=&k_{0}\int_{{\mbSigma}\in {\mathcal M}}~
     f({\mXbar}_{n}, {\bf S}_{n})|{\mbSigma}_{n}|^{-\frac{1}{2}(p+1)} ~~d{\mbSigma}   \\
   &=&k_{3}{|{\bf S}_{n}|}^{\frac{1}{2}(n-p-2)}{|{\bf V}_{n}|}^{-\frac{1}{2}n},  \nonumber
\end{eqnarray}
where the constant $k_{3}$ can be determinated by the equation
$\int_{{\cal D}} ~h_{2}({\mXbar}_{n},{\bf S}_{n})~~dx_{1}\cdots dx_{p}=1$.
Let $b_{2}(k,n,p)$ be defined the same as in (3.5) with $h_{2}({\mXbar}_{n},{\bf S}_{n})$
instead of $h_{1}({\mXbar}_{n},{\bf S}_{n})$. As such, we may obtain another critical points
for the LRT and the UIT. Further note that condition on ${\bf V}_{n}$, then (3.20) reduces
to that of Wang and McDermott (1998) for finding the critical value of semi-conditional LRT.

\indent From the above arguments, we may conclude that the phenomena whether the LRT and the
UIT are similar or not depend on the ways of defining the level of significance $\alpha$.
In passing, we note that different weight functions result in different critical values for
the LRT and the UIT respectively, as such the powers of the LRT and the UIT heavily depend on
the choice of weight functions for the nuisance parameter ${\mbSigma}$. Therefore, for a
specific test one may have several different power functions if different weight functions
are adopted. On the other hand, we may base on a specific weight function to study the power
properties and power dominances for some different tests. The phenomena of power
dominance for the interesting tests might be different under different weight functions.
The crux of using Bayesian notion to define the level of significance lies how to decide what
is the optimal weight function for the nuisance parameter. Thus, one of the disadvantages arises:
using Bayesian notions to define the level of significance may require extensive numerical
computations to get the critical values when $p\geq 4$.

\indent As a consequence, for a specific test when the rejected probability under null
hypothesis depends on the nuisance parameter ${\mbSigma}$, its power properties such as
similarity, unbiasedness and uniformly more powerful phenomenon relative to other tests may
be completely different if different ways of obtaining critical values are adopted.
In this paper, we regard these different ways of obtaining critical values for a specific test
as different kinds of definitions of significance level for that specific test. To give
a new definition of the level of significance can be viewed as a way to dig up the
insight information of null hypothesis. Of course, there are other approaches to define the
level of significance other than using Bayesian notions and maximum principle.

\indent In passing, by virtue of (3.5) and (3.6) we note that the Fisher approach for
$p$-values of the LRT and the UIT are hardly reported because their null distributions depend on
the nuisance parameter. For the problem of testing against the positive orthant space, the
main issue of reporting $p$-value still lies in how to deal with the nuisance parameter.
Several studies for handling the nuisance parameter in reporting $p$-value, which can also be
applied to the present problem, have been well investigated in the literature. We may refer to
Berger (1999) for the details. One of the main difficulties is how to find an unified way to
handle the problems. We are afraid that neither Fisher approach for reporting $p$-value alone
nor Neyman-Pearson approach alone is sound enough to take care the problem of testing
against the positive orthant space. Most people think that these two approaches disagree with
each other, however Lehmann (1993) asserted that these two approaches are complementary
rather than contradictory. In passing, we may note that the essential point of Fisher reporting
$p$-value is to find the optimal ways to obtain the insight information under null hypothesis, and
Neyman-Pearson theory essentially shows us how to find the tests that enjoy some optimal power
properties. As such, we may agree with Lehmann that there are essentially only one theory,
not two, for both Fisher and Neyman-Pearson approaches.

\vspace{0.3cm}
\def \theequation{4.\arabic{equation}}
\setcounter{equation}{0}

\indent {\bf 4. $D$-admissibility.} For the problem of testing against the positive orthant
space, the difficulty is that both the distribution functions of the LRT and the UIT under
null hypothesis are dependent on the nuisance paprameter. In the literature, many authors
heavily rely on Neyman-Pearson optimal theory to construct new tests so that their new tests
dominate the corresponding LRT in the sense of uniformly more powerful. For the details,
we refer to Wang and McDermott (1998), Berger (1989), Menendez and Salvador (1991),
Menendez et al. (1992) and the references therein. Sen and Tsai (1999) adopted Stein's
concept (1945) to establish the two-stage LRT, which enjoys some optimal properties,
however the cost is that the complicated distribution function of the test statistic arises. The method
proposed by Sen and Tsai (1999) is also essentially based on Neyman-Pearson optimal theory.
On the other hand, Perlman and Wu (1999) encouraged to abandon the criteria of
$\alpha$-admissibility, similarity and unbiasedness of Neyman-Pearson theory, and subjectively
suggested to use the LRT. Berger (1999) criticized Perlman and Wu's suggestion for not providing
any elegant support for the LRT. One of the main goals in this section is to seek a possible
resolution for the arguments between Perlman and Wu (1999) and Berger (1999).

\indent For testing against a global shift alternative the Hotelling's $T^{2}$-test is
uniformly most powerful invariant (UMPI), and hence, is also admissible (Simika, 1941). Stein
(1956) established the admissibility of the Hotelling's $T^{2}$-test by using the exponential
structure of the parameter space. The UMPI character, or even the admissibility of the
Hotelling's $T^{2}$-test may not generally hold for restricted alternatives, such as
$H_{1{\cal O}^{+}}$ in (1.4). The affine-invariance structure of the parameter space $\mathrm
{\Theta}=\{{\mbtheta}\in R^{p}\}$ does not hold for $H_{1{\cal O}^{+}}$, and hence, when
${\mbSigma}$ is arbitrary p.d. restriction to invariant tests makes little sense. As such, it
is conjectured, though not formally established, that possibly some other non-(affine)
invariant tests dominate Hotelling's $T^{2}$-tests, and hence, the latter is inadmissible.
Our contention is to establish that for testing $H_0$ vs $H_{1{\cal O}^{+}}$, the Hotelling's
$T^{2}$-test is inadmissible. For testing $H_0$ against $H_{1{\cal O}^{+}}$, the set of proper
Bayes tests and their weak limits only constitute a proper subset of essentially complete class
of tests (Marden, 1982). Marden's minimal complete classes consist of proper Bayes tests and
the tests with convex and decreasing acceptance regions for the density functions which are more general
than exponential family. For the problem (1.4), Theorem 2.1 tells us that the UIT is not a Bayes test.
Furthermore, the acceptance regions of the UIT may not be decreasing, though it can be shown
as the convex set. Hence, the conditions of Marden's minimal complete classes of tests are
also not satisfied for the UIT.

\indent We appraise Eaton's (1970) basic result on essentially complete class of test
functions for testing against restricted alternatives, when the underlying density belongs to
an exponential family as in our present context. Let $\Phi$ be Eaton's essentially complete
class of tests, so for any test ${\varphi}^{*}\notin {\Phi}$ there exists a test
${\varphi} \in {\Phi}$ such that ${\varphi}$ is at least as good as ${\varphi}^{*}$.

\vspace{0.3cm}
\indent {\sc Theorem 4.1.} {\it For the problem of testing $H_{0}: {\mbtheta} = {\bf 0}$ against
$H_{1{\cal O}^{+}}: {\mbtheta} \in {\cal O}_{p}^{+} \backslash \{{\bf 0}\}$,
the Hotelling's $T^{2}$-test is inadmissible.}

\indent {\sc Proof}. Note that the present testing hypothesis problem is invariant under the group of
transformations of positive diagonal matrices, hence, for simplicity, $\mbSigma$ may be treated
as the correlation matrix. Following Eaton (1970), we define
\begin{eqnarray}
 \Omega_{1}=\{{\mbSigma}^{-1}{\mbtheta}|~{\mbtheta} \in {\cal O}_{p}^{+}\}
 \backslash \{{\bf 0}\}.
\end{eqnarray}
Let ${\cal V}\subseteq R^{p}$ be the smallest closed convex cone containing $\Omega_{1}$.
Then the dual cone of ${\cal V}$ is defined as
\begin{eqnarray}
 {\cal V}^{-}=\{{\bf w}|<{\bf w},{\bf x}> \leq 0,~\forall~ {\bf x}\in {\cal V}\}.
\end{eqnarray}
Then, by (4.1) we have that ${\cal V}$ becomes a halfspace
and its dual cone is
\begin{eqnarray}
{\cal V}^{-}=\{{\mbtheta}\in R^{p}|~\sum_{i=1}^{p}{\theta}_{i}=0 \},
\end{eqnarray}
which is an unbounded hyperplane.

\indent The acceptance region of the Hotelling's $T^{2}$-test is given by
\begin{eqnarray}
\mathcal {A}_{T^{2}}=\{({\mXbar},{\bf S})|~T^{2} \leq {t^{2}_{\alpha}}\},
\end{eqnarray}
where ${t^{2}_{\alpha}}$ is the upper $100\alpha$\% point of the null hypothesis distribution
of $T^{2}$ (which is linked to a F-distribution). Since $\mathcal {A}_{T^{2}}$ is a
hyperellipsoidal set with origin ${\bf 0}$, it is bounded while $\mathcal {V}^{-}$,
as shown above, is unbounded. Therefore, Eaton's (1970) condition is not tenable. Hence
the Hotelling $T^{2}$-test is not a member of essentially complete class. ~~~~~~~~~~~~~Q.E.D.

\indent Next, for testing $H_0$ against $H_{1{\cal O}^{+}}$ we will show that the finite
union-intersection test (FUIT) (Tsai and Sen, 2004) and the union-intersection test (UIT)
(Sen and Tsai, 1999) are the members of Eaton's (1970) essentially complete class though they
are not the proper Bayes. A test which is a member of Eaton's (1970) essentially complete class
of tests can be established by showing the acceptance region contains the subspace
$\mathcal {V}^{-}$. Corresponding to a preassigned $\alpha~(0 < \alpha < 1)$, let $c_{\alpha}$
be the critical level, obtained by equating the right hand side of (3.9) to $\alpha$, thus the
UIT is a size-$\alpha$ test for $H_0$ vs $H_{1{\cal O}^{+}}$. Let ${\cal A}_{U}$ be the
acceptance region formed by letting $U_{n} \leq c_{\alpha}$ in (2.14). Partition the sample
space $R^{p}$ into $\cup_{a}I_{na}$, and for each $a$, {\reg}, let
${\cal A}_{a}=\{({\mXbar}_{n}, {\bf S}_{n})|n{\mXbar}'_{na:a'}
{\bf S}_{naa:a'}^ {-1}{\mXbar}_{na:a'} \leq c_{\alpha} \}I_{0na}$.
Treating ${\cal V}_{n}^{-}=\{({\mXbar}_{n}, {\bf S}_{n})|{\bf S}^{-1}_{n}{\mXbar}_{n}
\leq {\bf 0}\}=I_{0n\emptyset}$ as the skelton (pivotal set), then by (2.17)
we have that ${\cal A}_{U}={\cal V}_{n}^{-}\cup_{attach} {\cal A}_{a}$,
where $\cup_{attach}$ means that for each $a$, $\emptyset \subset a \subseteq P$,
the hyperspace ${\cal A}_{a}$ is attached to the boundary of ${\cal V}_{n}^{-}$ on the
subspace $I_{0na}$. By the property that ${\bf S}_{n}$ is positive definite with probability
one, we have ${\mathcal {V}^{-}} \subseteq \mathcal {V}^{-}_n \subseteq \mathcal {A}_{U}-{\bf a}^{0}$
for each ${\bf a}^{0}\in \partial {\mathcal {A}_{U}}$. Similarly, for the FUIT we may note
that its acceptance set ${\cal A}^{*}$ [see (2.21)] is a closed convex set and
$R^{-}_p =\{{\bf x}\in R^{p}| {\bf x} \leq {\bf 0}\}\subseteq {\cal A}^{*}$ as long as
$t_{n-1,{\alpha}^{*}}$ $\geq 0$ for ${\alpha}^{*} \leq \frac{1}{2}$.
Thus, the FUIT is a size-$\alpha$ test for $H_0$ vs $H_{1{\cal O}^{+}}$ and
${\mathcal {V}^{-}} \subseteq {\cal A}^{*}$.

\vspace{0.3cm}
\indent {\sc Theorem 4.2.}  {\it In the same setup of Theorem 4.1, both the UIT and the FUIT belong to Eaton's
essentially complete class of tests.}

\indent For the problem of testing $H_0$ against $H_{1{\cal O}^{+}}$, by virtues of Theorems 4.1 and 4.2 it is
interesting to see whether the Hotelling's $T^{2}$-test is power dominated by the UIT. For the case that both the
Hotelling's $T^{2}$-test and the UIT are unbiased, Tsai and Sen (2004) provided an affirmed answer for it when
$\mbSigma={\sigma}^{2} {\mbDelta}$, where ${\sigma}^{2}$ is an unknown scalar parameter and ${\mbDelta}$ is a known
$M$-matrix. However, the situation is still open when ${\mbSigma}$ is totally unknown. Theorem 4.2 does not guarantee
that the UIT (based on the test statistic ${U}_{n}$) is admissible. First, we examine the admissibility for the
corresponding UIT test based on the test statistic $U^{*}_{n}$ for the problem of testing against the halfspace.

\vspace{0.3cm}
\indent {\sc Theorem 4.3}. {\it For the problem of testing $H_{0}: {\mbtheta} = {\bf 0}$
against $H_{1 {\cal H}^{*}}: {\mbtheta} \in {\cal H}^{*}_{p} \backslash \{{\bf 0}\}$, the
UIT is $d$-admissible and $\alpha$-admissible.}

\indent {\sc Proof}. The acceptance region of UIT is
\begin{eqnarray}
  {\cal A}_{{U}^{*}}=\left\{({\mXbar}_{n}, {\bf S}_{n})|~ \sum_{\mbox{\rag}}
   \{ n{\mXbar}'_{na:a'}{\bf S}^{-1}_{naa:a'}{\mXbar}_{na:a'}\}
   I^{*}_{na}\leq k_{1}, ~{\bf S}_{n}~\hbox{positive definite}\right \},
\end{eqnarray}
where $P_{1}=P-\{p\}$ for a suitable $k_{1}$. Note that
\begin{eqnarray}
  {\cal A}_{{U}^{*}}={\cal A}_{U^{*}P} \cup {\cal A}_{U^{*}P1},
\end{eqnarray}
where
\begin{eqnarray}
  {\cal A}_{U^{*}P}=\{({\mXbar}_{n}, {\bf S}_{n})|~n{\mXbar}'_{n}{\bf S}^{-1}_{n}
  {\mXbar}_{n}\leq k_{1},~{\Bar X}_{np} \geq  0\},
\end{eqnarray}
and
\begin{eqnarray}
  {\cal A}_{U^{*}P1}=\{({\mXbar}_{n}, {\bf S}_{n})|~n{\mXbar}'_{nP_{1}: p}
   {\bf S}^{-1}_{nP_{1}P_{1}:p}{\mXbar}_{nP_{1}:p}
    \leq k_{1},~ \Bar{X}_{np} \leq 0 \}
\end{eqnarray}
respectively. For each $a, P_{1}\subseteq a \subseteq P$, Theorem 2.1 of Tsai (2003) tells us
that the function $\{ n{\mXbar}'_{na:a'}{\bf S}^{-1}_{naa:a'}{\mXbar}_{na:a'}\}$ is convex in
$({\mXbar}_{n},{\bf S}_{n})$. Thus the sets ${\cal A}_{U^{*}P}$ and ${\cal A}_{U^{*}P1}$ are
convex. Notice that the set ${\cal A}_{U^{*}P}$ is a half-ball in the space $\{({\mXbar}_{n},
{\bf S}_{n})|~{\Bar X}_{np} \geq  0\}$ and the set ${\cal A}_{U^{*}P1}$ is an unbounded
closed  cylinder-type convex cone in the space
$\{({\mXbar}_{n}, {\bf S}_{n})|~{\Bar X}_{np} \leq  0\}$. Also the extreme set (Rockafellar,
1972) of ${\cal A}_{U^{*}P}$ on the hyperplane
$\{({\mXbar}_{n}, {\bf S}_{n})|~{\Bar {X}_{np}}=0\}$ is identical to the one of
${\cal A}_{U^{*}P1}$ on the hyperplane $\{({\mXbar}_{n}, {\bf S}_{n})|~{\Bar {X}_{np}}=0\}$.
The sets ${\cal A}_{U^{*}P}$ and ${\cal A}_{U^{*}P1}$ are not disjoint since
${\cal A}_{U^{*}P1}={\cal A}_{U^{*}P}$ on the hyperplane
$\{({\mXbar}_{n}, {\bf S}_{n})|~{\Bar {X}_{np}}=0\}$, actually the intersection of these two
sets is identical to the extreme set of ${\cal A}_{U^{*}P}$ on the hyperplane
$\{({\mXbar}_{n}, {\bf S}_{n})|~{\Bar {X}_{np}}=0\}$. The set ${\cal A}_{{U}^{*}}$ is the
convex hull of ${\cal A}_{U^{*}P1}$ and ${\cal A}_{U^{*}P}$. Thus the set
${\cal A}^{*}_{U}$ is convex, the probability of $\partial {\cal A}_{{U}^{*}}$ is zero.

\indent Since ${\cal A}^{*}_{U}$ is convex on $({\mXbar}_{n}, {\bf S}_{n})$, we assume
${\cal A}^{*}_{U}$ is disjoint with the halfspace
\begin{eqnarray}
 {\mathcal H}_{c} = \{\, ({\mXbar}_{n}, {\bf S}_{n}) \mid {\mbomega}^{'}{\bf y}
    ={\mbnu}^{'}{\mXbar}_{n}-\frac{1}{2} tr{\mbLambda}{\bf V}_{n} > c, \}
\end{eqnarray}
where ${\mbLambda}$ is a symmetric matrix and ${\bf V}_{n}$ is p.d.. Let
${\cal H}^{-}_{p}=\{{\mbtheta}\in R^{p}|~\sum_{i=1}^{p}{\theta}_{i} \leq 0 \}$ and
${\cal H}^{\perp}_{p}$ be the halfspace which contains the hyperplane
${\cal H}^{0}_{p}=\{{\mbtheta}\in R^{p}|~\theta_{1}=\cdots=\theta_{p-1}=0, \theta_{p} < 0\}$
and is perpendicular to the hyperplane
${\cal H}^{1}_{p}=\{{\mbtheta}\in R^{p}|~\sum_{i=1}^{p}{\theta}_{i}=0 \}$. For $\mbSigma$
being any arbitrarily positive definite matrices,
${\mbSigma}^{-1}: R^{p(p+1)/2} \to R^{p(p+1)/2}$ is a continuous, linear and open
mapping on the general linear group $GL(R^{p(p+1)})$. Let
${\cal S}={\cal H}^{-}_{p} \cap {\cal H}^{\perp}_{p}$, thus
${\mbtheta} \in {\cal H}^{*}_{p}$ implies that
${\mbnu}=n{\mbSigma}^{-1}{\mbtheta} \in \overline {R^{p}\backslash{\cal S}}$.
Here we want to show that $\mbomega_{0}+{\lambda}{\mbomega} \in (\overline {R^{p}
\backslash{\cal S}})\backslash \{{\bf 0}\}$ holds, namely to show
${\bf I}+{\lambda}{\mbLambda}$ is positive definite (it suffices to show that
${\mbLambda}$ is positive semidefinite) and
$\mbnu_{0}+{\lambda} {\mbnu} \in (\overline {R^{p}\backslash{\cal S}})\backslash \{{\bf 0}\}$
for ${\lambda} > 0$.

\indent (I). To show that ${\mbLambda}$ is positive semidefinite. Suppose that ${\mbLambda}$
is not, and then assume that it is
\[
{\mbLambda}= {\bf D}\left[\begin{array}{cccc}
    {\bf I}_{b_{1}}   &  {\bf 0 }          &  {\bf 0}   \\
    {\bf 0}   & -{\bf I}_{b_{2}}   &  {\bf 0 }  \\
    {\bf 0}   &  {\bf 0}           &  {\bf 0}_{b_{3}}
\end{array} \right]{\bf D}^{'}_,
\]
where ${\bf D}$ is nonsingular, and $1 \leq b_{2} < p, 1 \leq b_{1}+b_{2} <p$ and
$\sum_{i=1}^{3}b_{i}=p$. Take ${\mXbar}_{n}=\frac{1}{{\lambda}}{\mZbar}_{n}$ and
\begin{eqnarray}
  {\bf V}_n = {\lambda}({\bf D}^{'})^{-1}\left[\begin{array}{cccc}
    {\bf I}_{b_{1}}   &  {\bf 0 }          &  {\bf 0}   \\
    {\bf 0}   & {\bf I}_{b_{2}}   &  {\bf 0 }  \\
    {\bf 0}   &  {\bf 0}           &  {\bf I}_{b_{3}}
\end{array} \right]{\bf D}^{-1}
\end{eqnarray}
for positive $\lambda$. Then
\[
{\mbomega}^{'}{\bf y}=\frac{1}{{\lambda}}{\mbnu}^{'}{\mZbar}_{n}-\frac{{\lambda}}{2}~tr
\left[\begin{array}{cccc}
    {\bf I}_{b_{1}}   & {\bf 0 }          &{\bf 0}   \\
    {\bf 0}   & -{\bf I}_{b_{2}}    &{\bf 0 }  \\
    {\bf 0}   & {\bf 0}          &{\bf 0}_{b_{3}}
\end{array} \right], \]
which is greater than $c$, a positive constant. Let the space
$\mathcal D^{*}=\{({\mXbar}_{n}, {\bf S}_{n})|~{\bf S}_{n}~\mbox {is p.d.}~
\linebreak \mbox {and}~{\bf V}_n~\mbox {has the form in (4.10)}\}$. Thus, the halfspace
in (4.9) becomes
\begin{eqnarray}
\mathcal {H}_{c} \cap \mathcal D^{*}=
   \{{\mZbar}_{n}|~{\mbnu}^{'}{\mZbar}_{n} > {\lambda} c
   -\frac {{\lambda}^{2}}{2}(b_{1}-b_{2})\}.
\end{eqnarray}
 Let
\begin{eqnarray}
{\cal A}_{T^{2}, k_{1}}=\{({\mXbar}_{n}, {\bf S}_{n})|~n{\mXbar}'_{n}{\bf S}^{-1}_{n}
  {\mXbar}_{n}\leq k_{1}\}.
\end{eqnarray}
Note that
\begin{eqnarray}
  {\cal A}_{{U}^{*}} \supseteq {\cal A}_{T^{2}, k_{1}}
   =\{({\mXbar}_{n}, {\bf S}_{n})|~n{\mXbar}'_{n}{\bf V}^{-1}_{n}{\mXbar}_{n}\leq k^{*}_{1}\}
\end{eqnarray}
for a suitable positive constant $k_{2}$. Thus,
\begin{eqnarray}
{\cal A}_{T^{2}, k_{1}} \cap \mathcal D^{*}=\{{\mZbar}_{n}|~n{\mZbar}'_{n}
{\bf D}{\bf D}^{'}{\mZbar}_{n}\leq {\lambda}^{3}k^{*}_{1}\}.
\end{eqnarray}
By virture of (4.11), (4.13) and (4.14), the intersection of
${\cal A}_{{U}^{*}} \cap \mathcal D^{*}$
and $\mathcal {H}_{c}\cap \mathcal D^{*}$ is nonempty for sufficiently large
${\lambda}$. This contradicts the fact that (4.5) and (4.9) are disjoint. Thus
${\mbLambda}$ is positive semidefinite.

\indent (II). To show that $\mbnu_{0}+{\lambda} {\mbnu} \in (\overline {R^{p}
\backslash{\cal S}})\backslash \{{\bf 0}\}$ for ${\lambda} > 0$. Our aim is to show that
${\mbnu} \in (\overline {R^{p}\backslash{\cal S}})\backslash \{{\bf 0}\}$. To proceed, we
shall work with another halfspace
\begin{eqnarray}
{\mathcal H}^{a}_{c} = \{\, ({\mXbar}_{n}, {\bf S}_{n}) \mid n\mbnu'{\mXbar}_{n}
     -\frac{1}{2} \tr [({\bf I}+{\mbLambda})({\bf S}_{n}
    + n{\mXbar}_{n}{\mXbar}^{'}_{n})] > c \,\}.
\end{eqnarray}
Note that ${\bf S}_{n}$ is p.d. with probability one, thus
$(\mbGamma, {\bf S}_{n}) \in {\mathcal H}^{a}_{c}$ implies that
$({\mXbar}_{n}, {\bf S}_{n}) \in {\mathcal H}_{c}$. Hence, we have
\begin{eqnarray}
{\mathcal H}^{a}_{c} \subset {\mathcal H}_{c}.
\end{eqnarray}
Consequently, the intersection of ${\mathcal A}_{{U}^{*}}$ and ${\mathcal H}^{a}_{c}$ is
empty too. Note that ${\bf I}+{\mbLambda}$ is p.d., thus we may denote its inverse by
${\mbGamma}$, i.e., $({\bf I}+{\mbLambda})^{-1}={\mbGamma}$. Let
$\Tilde{\mbnu}={\mbGamma}{\mbnu}$ and ${\mZbar}_{n}={\mbGamma}^{-1}{\mXbar}_{n}$, then
${\mathcal H}^{a}_{c}$ becomes
\begin{eqnarray}
{\mathcal H}^{a}_{c} = \{\, ({\mZbar}_{n}, {\bf S}_{n}) \mid
  n \Tilde{\mbnu}'{\mZbar}_{n}-\frac{1}{2} n{\mZbar}^{'}_{n}{\mbGamma}{\mZbar}_{n}
  -\frac{1}{2} \tr ({\mbGamma}^{-1}{\bf S}_{n}) > c \,\}.
\end{eqnarray}
And then ${\mbnu} \in (\overline {R^{p}\backslash{\cal S}})\backslash \{{\bf 0}\}$ implies
that $\Tilde{\mbnu} \in R^{p}\backslash \{{\bf 0}\}$, hence our aim becomes to show that
$\Tilde{\mbnu} \in R^{p}\backslash \{{\bf 0}\}$. Suppose that the statement is not true,
then we have to show that $\Tilde{\mbnu}={\bf 0}$ leads to a contrdiction.

\indent Take ${\bf S}_{n}={\mbGamma}$ and let $\Tilde{\mbnu}={\bf 0}$, then we have
\begin{eqnarray}
{\mathcal H}^{a}_{c}=\{{\mZbar}_{n} \mid n{\mZbar}^{'}_{n}{\mbGamma}{\mZbar}_{n}
   \leq  -(p+2c) \}.
\end{eqnarray}
Also note that given ${\bf S}_{n}={\mbGamma}$, ${\cal A}_{T^{2}, k_{1}}$ becomes
\begin{eqnarray}
  {\cal A}_{T^{2}, k_{1}}
   =\{{\mZbar}_{n} \mid n{\mZbar}^{'}_{n}{\mbGamma}{\mZbar}_{n}\leq 2k^{*}_{1}\}.
\end{eqnarray}
It is obvious that ${\mathcal H}^{a}_{c}\cap {\cal A}_{T^{2}, k_{1}} \ne \emptyset$, and hence
${\mathcal H}^{a}_{c}\cap {\cal A}_{U^{*}} \ne \emptyset$. This contradicts the fact that (4.5) and (4.9) are disjoint.
Hence the condition that $\mbomega_{1}+{\lambda}{\mbomega}\in {\cal H}^{*}_{p}$ for an arbitrarily large ${\lambda}$
is satisfied, namely the conditions of Theorem 8 of Lehmann (1986) hold. Thus by Theorem 8 of Lehmann (1986),
we conclude that for testing against the halfspace the UIT is $d$-admissible.

\indent With similar arguments as in Theorems 2.1 and 3.1 of Sen and Tsai (1999), we may show that the UIT is similar
and unbiased. Therefore by Corollary 2 of Lehmann (1986), the UIT is $\alpha$-admissible. ~~~~~~~~~~~~~Q.E.D.

\indent For the problem of testing $H_0$ against $H_{1{\cal O}^{+}}$, it is interesting to see whether the UIT are
admissible in decision-theoretic sense. Our problems are under the set up of exponential family, Theorem 8 of
Lehmann (1986), which states the conditions for tests to be $d$-admissible for the problems of testing against
restricted alternatives, can be viewed as a generalization of Birnbaum (1955) and Stein (1956). Based on the results
of Theorem 8, Corollary 2 of Lehmann (1986) and Theorem 2.1 of Tsai (2003) we prove the $d$-admissibilities of UIT
for the problem (1.4). With arguments parallel to those in the proof of Theorem 4.3, we also have the following.

\vspace{0.3cm}
\indent {\sc Theorem 4.4}. {\it For the problem of testing $H_{0}: {\mbtheta} = {\bf 0}$
against $H_{1{\cal O}^{+}}: {\mbtheta} \in {\cal O}_{p}^{+} \backslash \{{\bf 0}\}$, the UIT
is $d$-admissible.}

\indent {\sc Outline of the proof}. The acceptance region of UIT is
\begin{eqnarray}
  {\cal A}_{U}=\left\{({\mXbar}_{n}, {\bf S}_{n})|~\sum_{\mbox{\reg}}\{n{\mXbar}'_{na:a'}
  {\bf S}_{naa:a'}^{-1}{\mXbar}_{na:a'}\}I_{na}
    \leq k_{2},~{\bf S}_{n}~\hbox{positive definite} \right \},
\end{eqnarray}
for a suitable $k_{2}$. Note that
\begin{eqnarray}
  {\cal A}_{{U}}=\cup_{{\mbox{\reg}}}{\cal A}_{Ua},
\end{eqnarray}
where
\begin{eqnarray}
  {\cal A}_{Ua}=\{({\mXbar}_{n}, {\bf S}_{n})|~n{\mXbar}'_{na:a'}
    {\bf S}_{naa:a'}^{-1}{\mXbar}_{na:a'}\leq k_{2},
   ~{\mXbar}_{na:a'}>{\bf 0}, {\bf S}_{na'a'}^{-1}{\mXbar}_{na'}\leq 0\}.
\end{eqnarray}
For each $a,~{\mbox{\reg}} $, the function
$\{ n{\mXbar}'_{na:a'}{\bf S}^{-1}_{naa:a'}{\mXbar}_{na:a'}\}I_{na}$ is convex in
$({\mXbar}_{n},{\bf S}_{n})$ for ${\bf S}_{n}$ positive definite. Thus each ${\cal A}_{Ua}$
is convex. Notice that any two sets ${\cal A}_{Ua}$ and ${\cal A}_{Ub}$ are not disjoint,
where $\emptyset \subseteq a,~ b \subseteq P$ with $||a|-|b||=1$ and the intersection of
these two sets is identical to their common extreme set. Treating
${\cal A}_{U\emptyset}=\{({\mXbar}_{n}, {\bf S}_{n})|~{\bf S}_{n}^{-1}{\mXbar}_{n}\leq {\bf 0}\}$
as the skeleton, we may note that
${\cal A}_{U}={\cal A}_{U\emptyset}\cup_{attach} {\cal A}_{Ua}$, where $\cup_{attach}$
means that for each $a$, $\emptyset \subset a \subseteq P$,
the hyperspace ${\cal A}_{Ua}$ is attached to the common extreme set with the set of ${\cal A}_{U\emptyset}$.
The set ${\cal A}_{{U}}$ is the convex hull of ${\cal A}_{Ua},~\emptyset \subset a \subseteq P$. Thus the set
${\cal A}_{{U}}$ is convex; the probability of $\partial {\cal A}_{{U}}$ is zero.

\indent Assume ${\cal A}_{{U}}$ is disjoint with the halfspace in (4.9). Note that
${\mbtheta}\in {\cal O}^{+}_{p}$ implies that
${\mbnu}=n{\mbSigma}^{-1}{\mbtheta} \in \overline{R^{p}\backslash {\cal H}^{-}_{p}}$, where
${\cal H}^{-}_{p}=\{{\mbtheta}\in R^{p}|~\sum_{i=1}^{p}{\theta}_{i} \leq 0 \}$. Here we want
to show that
$\mbomega_{1}+{\lambda}{\mbomega} \in (\overline{R^{p}\backslash {\cal H}^{-}_{p}}) \backslash \{{\bf 0}\}$ holds,
i.e., $\mbnu_{1}+{\lambda} {\mbnu} \in (\overline{R^{p}\backslash {\cal H}^{-}_{p}})
\backslash \{{\bf 0}\}$ and ${\mbLambda}$ is positive semidefinite for $\lambda > 0$.
We shall take ${\mbnu}_{1}={\bf 0}$. With similar arguments as in (1) of Theorem 4.3, we can
claim that ${\mbLambda}$ is positive semidefinite. We also adopt the notations
$\Tilde {\mbnu}$, ${\mZbar}_{n}$ and ${\mathcal H}^{a}_{c}$ from Theorem 4.3. Thus,
${\mbnu} \in (\overline{R^{p}\backslash {\cal H}^{-}_{p}})$ implies that
$\Tilde {\mbnu} \in \overline{R^{p}\backslash {\cal R}^{-}_{p}}$, where
$R^{-}_{p}=\{{\bf x} \in R^{p} \mid {\bf x} \leq {\bf 0}\}$.

\indent Take ${\bf S}_{n}={\mbGamma}$, then we have
\begin{eqnarray}
{\mathcal H}^{a}_{c}&=&\{{\mZbar}_{n} \mid  2n \Tilde{\mbnu}'{\mZbar}_{n}
          -n{\mZbar}^{'}_{n}{\mbGamma}{\mZbar}_{n} \geq  p+2c \}   \\
   &=&  \{{\mZbar}_{n} \mid n({\mZbar}_{n}-{\mbGamma}^{-1}\Tilde{\mbnu})^{'}{\mbGamma}
    ({\mZbar}_{n}-{\mbGamma}^{-1} \Tilde{\mbnu}) \leq
     n {\Tilde{\mbnu}}^{'}{\mbGamma}^{-1} {\Tilde{\mbnu}}-(p+2c) \} \nonumber
\end{eqnarray}
Our aim is to show that
$\Tilde{\mbnu} \in (\overline{R^{p}\backslash {\cal R}^{-}_{p}})\backslash \{{\bf 0}\}$.
Suppose it is not true, then $\Tilde{\mbnu} \in {R}^{-}_{p}$ and hence
${\mathcal H}^{a}_{c}$ contains an open ball in the space
$\{{\mZbar}_{n}\mid {\mZbar}_{n} \leq {\bf 0}\}$. Given ${\bf S}_{n}={\mbGamma}$, it is easy
to see that ${\cal A}_{U\emptyset}=\{{\mZbar}_{n}\mid {\mZbar}_{n}\leq {\bf 0}\}$. Thus, we
have ${\mathcal H}^{a}_{c}\cap {\cal A}_{T^{2}, k_{1}} \ne \emptyset$, and
hence ${\mathcal H}^{a}_{c}\cap {\cal A}_{U^{*}} \ne \emptyset$. This contradicts the
fact that (4.5) and (4.9) are disjoint. Hence the condition that
$\mbomega_{1}+{\lambda}{\mbomega} \in (\overline{R^{p}\backslash {\cal H}^{-}_{p}})
\backslash \{{\bf 0}\}$ for an arbitrarily large ${\lambda}$ is satisfied, namely the
conditions of Theorem 8 of Lehmann (1986) hold. Therefore by Theorem 8 of Lehmann (1986),
we conclude that for testing against the positive orthant space the UIT is $d$-admissible. ~~~~~~~~~~~~~Q.E.D.

\indent For the problem of testing $H_0$ against $H_{1{\cal O}^{+}}$, Proposition 3.1 shows that the UIT based on
the test statistic $U_{n}$ is $\alpha$-inadmissible because its power is dominated by that of the corresponding UIT
based on the test statistic $U^{*}_{n}$, this result is against our common statistical sense. This is due to the fact
that the power function of UIT based on $U_{n}$ under null hypothesis depends on the nuisance parameter${\mbSigma}$
as discussed in above section.  To avoid this unpleasant phenomeno in this situation,  the spirits lies in the
compromise of Fisher approach emphasing on the type I erros and Neyman-Pearson optimal theory without detailed
consideration of power and power comparisions of tests would make sense only in the situation when the tests are
similar. As such, Wald's decision theory via $d$-admissibility and Neyman-Pearson's optimal theory is essentially one
theory, not two. This way can also provide a resolution for the arguments between Perlman and Wu (1999) and
Berger (1999). Therefore for the problem of testing $H_0$ against $H_{1{\cal O}^{+}}$, we suggest use the unified
conservative maximum principle to define the level of significance and then adopt the criterion in decision-theoretic
sense (d-admissibility) for the power dominance problems.

\vspace{0.3cm}
\def \theequation{5.\arabic{equation}}
\setcounter{equation}{0}

\indent {\bf 5. The difficulty of LRT}. When ${\mbSigma}$ is totally unknown, Proposition 2.1
fails to characterize whether or not the LRT is the generalized Bayes test. Also, let $l_{\alpha}$
be the critical point of the level of significance $\alpha$ for the LRT for the problem of
testing against the positive orthant space, which can be obtained from the right hand sides of
(3.8), and ${\mathcal A}_{L}=\cup_{\emptyset}^{P}{\cal A}_{La}$ be the acceptance region of the
LRT, where ${\cal A}_{La}=\{({\mXbar}_{n}, {\bf S}_{n})| ~n{\mXbar}'_{na:a'} {\bf S}^{-1}_{naa:a'}
{\mXbar}_{na:a'}/(1+n {\mXbar}'_{na'}{\bf S}^{-1}_{na'a'}{\mXbar}_{na'})\leq l_{\alpha}\}
I_{na}$, {\reg}. For simplicity, we take $p=2$ and ${\bf S}_{n}$ be the diagonal matrix, then
it is easy to see that the acceptance region ${\mathcal A}_{L}$ is not convex. Similar
arguments, it is true that the acceptance region of the LRT for testing against the
half-space is not convex. In the literature, most of the results of admissibility for tests are
established based on the fundamental assumption that the acceptance region is convex. The
nonconvexity of acceptance region of the LRT may cause the difficulty for the question whether
the LRT is $d$-admissible for problem (1.4). We will put it as the project for future study.

\vspace{0.3cm}

\vspace {0.3cm}
\centerline{\ REFERENCES}
\begin{description}

\item {\sc Anderson, T. W.} (2003). {\it An Introduction to Multivariate Statistical
Analysis}, 3rd edition. Wiely, New York.

\item {\sc Basu, D.} (1977). On the elimination of nuisance parameters. {\it J.
Amer. Statist. Assoc.} {\bf 72} 355-366.

\item {\sc Berger, J. O., Liseo, B.} and {\sc Wolpert, R. L.} (1999). Integrated likelihood
methods for eliminating nuisance parameters. {\it Statist. Sci.} {\bf 14} 1-28.

\item {\sc Berger, R. L.} (1989). Uniformly more powerful tests for hypotheses concerning
linear inequalities and normal means. {\it J. Amer. Statist. Assoc.} {\bf 84} 192-199.

\item {\sc Berger, R. L.} (1999). Comment on ``The emperor's new tests", {\it Statist.
Sci.} {\bf 14} 370-373.

\item {\sc Birnbaum, A.} (1955). Characterization of complete classes of tests of some
multiparametric hypothesis, with applications to likelihood ratio tests. {\it Ann. Math.
Statist.} {\bf 26} 21-36.

\item {\sc Brown, L. D.} (1986). {\it Fundamentals of Statistical Exponential families,
with Application in Statistical Decision Theory}. Institute of Mathematical Statistics,
Lecture Notes-Monograph Series, Volume {\bf 9}. Hayward, california.

\item {\sc Chang, T.} and {\sc Eaves, D.} (1990). Reference priors for the orbit in a group
model. {\it Ann. Statist.} {\bf 18} 1595-1614.

\item {\sc Eaton, M. L.} (1970). A complete class theorem for multidimensional one sided
alternatives. {\it Ann. Math. Statist.} {\bf 41} 1884-1888.

\item {\sc Kiefer, J.} and {\sc Schwartz, R.} (1965). Admissible Bayes character of $T^{2}-$,
$R^{2}-$, and other fully invariant tests for classical multivariate normal problems.
{\it Ann. Math. Statist.} {\bf 36} 747-770.

\item {\sc Lehmann, E. L.} (1986). {\it Testing Statistical Hypotheses}, 2nd edition.
Wiely, New York.

\item {\sc Lehmann, E. L.} (1993). The Fisher, Neyman-Pearson theories of testing
hypotheses: One theory or two? {\it J. Amer. Statist. Assoc.} {\bf 88} 1242-1249.

\item {\sc Marden, J. I.} (1982). Minimal complete classes of tests of hypothesis with
multivariate one-sided alternatives. {\it Ann. Statist.} {\bf 10} 962-970.

\item {\sc Menendez, J. A.} and {\sc Salvador, B. } (1991). Anomalies of the likelihood
ratio test for testing restricted hypotheses. {\it Ann. Statist.} {\bf 19} 889-898.

\item{\sc Menendez, J. A.}, {\sc Rueda, C.} and {\sc Salvador, B.} (1992). Dominance
of likelihood ratio tests under cone constraints. {\it Ann. Statist.} {\bf 20} 2087-2099.

\item {\sc N\"{u}esch, P. E.} (1966). On the problem of testing location in multivariate
problems for restricted alternatives. {\it Ann. Math. Statist.} {\bf 37} 113-119.

\item {\sc Perlman, M. D.} (1969). One-sided problems in multivariate analysis. {\it
Ann. Math. Statist.} {\bf 40} 549-567.

\item {\sc Perlman, M. D.} and {\sc Wu, L.} (1999). The emperor's new tests. {\it Statist.
Sci.} {\bf 14} 355-369.

\item {\sc Rockafellar, R. T.} (1972). {\it Convex Analysis}. Princeton University Press.

\item {\sc Roy, S. N.} (1953). On a heuristic method of test construction and its use in
multivariate analysis. {\it Ann. Math. Statist.} {\bf 24} 220-238.

\item {\sc Roy, S. N., Gnanadesikan, R.} and {\sc Srivastava, J.N.} (1972). {\it Analysis
and Deign of Certain Quantitative Multiresponse Experiments}. Pergamon Press, Oxford.

\item {\sc Sen, P. K.} and {\sc Tsai, M. T.} (1999). Two-stage LR and UI tests for one-sided
alternatives: multivariate mean with nuisance dispersion matrix. {\it J. Multivariate Anal.}
{\bf 68} 264-282.

\item {\sc Simaika, J. B.} (1941). On an optimum property of two important statistical
tests. {\it Biometrika}  {\bf 32} 70-80.

\item {\sc Stein, C.} (1945). A two-sample test for a linear hypothesis whose power function
is independent of $\sigma$. {\it Ann. Math. Statist.} {\bf 16} 243-258.

\item {\sc Stein, C.} (1956). The admissibility of Hotelling's $T^{2}$-test. {\it Ann. Math.
Statist.} {\bf 27} 616-623.

\item {\sc Tang, D. I.} (1994). Uniformly more powerful tests in a one-sided multivariate
problem. {\it J. Amer. Statist. Assoc.} {\bf 89} 1006-1011.

\item {\sc Tsai, M. T.} (2003). On the invariant tests for means with covariates. {\it
TR 2003-05}. Technical Report, Institute of Statistical Science, Academia Sinica.

\item {\sc Tsai, M. T.} and {\sc Sen, P. K.} (2004). On inadmissibility of Hotelling's
$T^{2}$-test for restricted alternatives. {\it J. Multivariate Anal.} {\bf 89} 87-96.

\item {\sc Wald, A.} (1950). {\it Statistical Decision Functions}. Wiely, New York.

\item {\sc Wang, Y.} and {\sc McDermott, M. P.} (1998). Conditional likelihood ratio test
for a nonnegative normal mean vector. {\it J. Amer. Statist. Assoc.} {\bf 93} 380-386.

\end{description}

\vspace{0.3cm}
\newlength{\lewidth}
\newlength{\riwidth}
\settowidth{\lewidth}{University of Virginia at Charlottesville, Virginia}
\settowidth{\riwidth}{Institute of Statistical Science at Taipei}

\begin{minipage}[t]{\riwidth}
{\sc Institute of Statistical Science }\\
{\sc Academia Sinica}  \\
{\sc Taipei  11529, Taiwan}   \\
{ E-mail: mttsai@stat.sinica.edu.tw}   \\
\end{minipage}

\end{document}